\documentclass[10pt]{article}
\usepackage[utf8]{inputenc}
\usepackage[a4paper,left=2cm,right=1.5cm, top=2.5cm,bottom=2.5cm]{geometry}
\usepackage{pgfplots}
\usepackage{subcaption}
\usepackage{adjustbox}
\usepackage{graphicx, tikz-3dplot}
\usepackage{pgfplotstable}
\usepackage[justification=centering]{caption}
\usepackage{bm}
\usepackage{dsfont}

\usepackage{mathtools}
\usepackage{amsmath}
\usepackage{amssymb}

\usepackage{multicol}
\usepackage{authblk}

\usepackage{comment}
\usepackage[numbers,super]{natbib}

\usepackage{algorithm,algpseudocode}
\usepackage{hyperref}
\newtheorem{remark}{Remark}

\definecolor{EDFblue}{RGB}{5,58,115}
\definecolor{EDForange}{RGB}{236,78,4}
\definecolor{EDFgreen}{RGB}{0,140, 33}

\newcommand{\norm}[1]{\left\lVert#1\right\rVert}

\title{A projection-based reduced-order model for parametric quasi-static nonlinear mechanics using an open-source industrial code}

\author[1,2,3]{Eki Agouzal}
\author[1]{Jean-Philippe Argaud}
\author[2,3]{Michel Bergmann}
\author[1]{Guilhem Ferté}
\author[2,3]{Tommaso Taddei}
\affil[1]{EDF Lab Paris-Saclay, EDF R\&D, 7 Boulevard Gaspard Monge, 91120 Palaiseau, France}
\affil[2]{IMB, UMR 5251, Univ. Bordeaux, 33400 Talence, France}
\affil[3]{INRIA, Inria Bordeaux Sud-Ouest, Team MEMPHIS, Univ. Bordeaux, 33400 Talence, France}
\date{December 2022}

\usepackage{fancyhdr}
\fancypagestyle{firstpage}
{
    \fancyhead[L]{}    
    \fancyhead[R]{Submitted to International Journal for Numerical Methods in Engineering, December 2022}
}

\begin{document}

\maketitle
\thispagestyle{firstpage}

\begin{abstract}
We propose a projection-based model order reduction procedure for a general class of parametric quasi-static problems in nonlinear mechanics with internal variables. The methodology is integrated in the industrial finite element code \textsf{code$\_$aster}. Model order reduction aims to lower the computational cost of engineering
studies that involve the simulation to a costly high-fidelity differential model for many different parameters, which correspond, for example to material properties or initial and boundary conditions. We develop an adaptive algorithm based on a POD-Greedy strategy, and we develop an hyper-reduction strategy based on an element-wise empirical quadrature in order to speed up the assembly costs of the reduced-order model by building an appropriate reduced mesh.  We introduce a cost-efficient error indicator which relies on the reconstruction of the stress field by a Gappy-POD strategy. We present numerical results for a three-dimensional elastoplastic system in order to illustrate and validate the methodology.
\end{abstract}

\section{Introduction}\label{sec1}

\subsection{Context}

Numerical simulations have been used for a long time within engineering studies, often in the perspective of evaluating the same study for slightly different configurations. These variations may include changes in the input signals, in the actual model parameters, or even in the geometry (\emph{many-query} problem for parametric studies). For problems modeled by partial differential equations (PDEs), extensive explorations of the parameter domain based on standard finite element (FE) solvers are prohibitively expensive. Model order reduction (MOR\cite{hesthaven2016certified, quarteroni2015reduced, rozza2008reduced}) consists in a broad spectrum of algorithms that aim to drastically reduce the marginal cost associated with one computation, by taking into account prior knowledge from previous high-fidelity simulations. Parametric model order reduction (pMOR) refers to a class of techniques that aim at constructing a low-dimensional surrogate (or reduced-order) model (ROM) to approximate the solution field over a range of parameters, by taking into account prior knowledge from previous high-fidelity (HF) simulations.

Our aim is to devise an intrusive pMOR procedure for large-scale problems in non-
linear structural mechanics that is consistent with an industrial code used in practice by engineers for HF simulations. Intrusive pMOR techniques rely on the projection of the differential operator onto suitable empirical reduced spaces, and thus require the access to local assembly routines of the underlying HF code. Intrusive techniques need to be elaborated in compliance with the operators and data structures used in the HF industrial code: the key challenge is to benefit from the robustness of the pre-existing industrial code — which allows to run real-world simulations for three-dimensional complex geometries and non-trivial mechanical behaviors — without having to modify the overall architecture (i.e., data structures and local assembly routines) of the HF code. In this work, we focus on the open-source software \textsf{code$\_$aster}\cite{aster}: \textsf{code$\_$aster} is a well established, qualified and broadly-used industrial grade finite element solver for structural mechanics studies that is mainly developed within Electricit{\'e} De France (EDF)’s R\&D.

In this contribution, we focus on a general class of parametric mechanical problems with internal variables in a nonlinear quasi-static framework, where we consider small-displacement small-strain mechanical problems.

\subsection{Objective of the paper and relation to previous works}
\label{sec:intro:subsec:purpose}

The main contribution of this work is the formulation and implementation of an hyper-reduced model for nonlinear quasi-static mechanical problems based on the industrial finite element code \textsf{code$\_$aster}\cite{aster}. We develop an adaptive algorithm, whose design is rooted in the offline-online paradigm. The algorithm can be divided in two steps : an \emph{offline} (or training) step, where a set of basis function is built from  a database of several HF solutions in order to approximate the solution manifold, and an \emph{online} step, during which the approximate solution is sought for a new set of parameter values. The algorithm we hereby present is founded on a Proper Orthogonal Decomposition (POD)-Greedy strategy, which was introduced in Reference \citenum{haasdonk2008reduced} and analyzed in Reference \citenum{haasdonk2013convergence}. Similarly to the weak-Greedy algorithm for stationary problems\cite{buffa2012priori, binev2011convergence}, the POD-Greedy procedure iteratively explores the parameter domain to identify poorly-approximated configurations through the vehicle of an a posteriori error indicator, and relies on the Proper Orthogonal Decomposition (POD\cite{berkooz1993proper, bergmann2009enablers, volkwein2011model}) to compress the temporal trajectory.  In this work, we rely on a time-averaged error indicator in a similar way to what has been done in Reference \citenum{iollo2022adaptive}, inspired by Reference \citenum {fick2018stabilized}.

Our solution strategy relies on a Galerkin projection method. Since the operator is nonlinear, the computational complexity of the operator assembly (jacobian and residuals) scales with the size of the HF model. In order to circumvent this obstacle, we develop an hyper-reduction strategy based on empirical quadrature (EQ) : our approach relies on the construction of a reduced mesh to speed up online assembly costs of the ROM. We refer here to a reduced mesh to describe a mesh designed by considering a subset of the cells of the HF mesh. The EQ procedure has been first proposed in References \citenum{farhat2014dimensional, farhat2015structure} and used in several previous work\cite{riffaud2021dgdd, iollo2022adaptive}. This approach relies on the reweighting of either the quadrature points of the mesh\cite{yano2019lp}, or the elemental contributions\cite{yano2019discontinuous, iollo2022adaptive}, in order to approximate the residuals. Several other techniques have been introduced in the literature in order to dodge the bottleneck induced by the projection step for nonlinear non-affine problems. Other reweighting methods have been introduced such as the Empirical Cubature Method\cite{hernandez2017dimensional}, which inspired implementation within industrial context\cite{casenave2020nonintrusive, casenave2019error}. Hyper-reduction approaches also include the family of algorithms derived from the Empirical Interpolation Method\cite{barrault2004empirical}, which encompass its discrete variant\cite{chaturantabut2010nonlinear}, or techniques which belong to the Gappy-POD application, such as the A priori Hyper-Reduction\cite{ryckelynck2005priori}, or the Gauss Newton with approximated tensors\cite{carlberg2013gnat}.

During the past decade, advances in MOR have led to the application of online-efficient projection-based ROMs to a broad range of problems in mechanics. In more details, several authors have considered the application of projection-based MOR techniques to large-scale three-dimensional problems in nonlinear mechanics including contact\cite{ballani2018component}\cite{le2022condition} , thermo-mechanics\cite{lindsay2022preconditioned} , and elasto-viscoplasticity\cite{casenave2020nonintrusive}\cite{casenave2019error}.

Our work is a continuation of the research effort carried out at EDF R\&D to deploy effective ROMs for nonlinear problems in structural mechanics. In this respect, we mention earlier works on nonlinear parabolic thermo-mechanical problems\cite{benaceur2018reduction} , on vibro-acoustics problems\cite{khoun2021reduced}, and also on welding \cite{dinh2018modeles} ; in particular, the work in Reference \citenum{dinh2018modeles} represents one of the first efforts to devise hyper-reduced ROMs in \textsf{code$\_$aster}. As discussed in section \ref{sec:formulation}, the \textsf{code$\_$aster} framework involves dualization of the boundary conditions and relies on a mesh hierarchy that comprises a three-dimensional mesh — for volumetric terms — and a two-dimensional mesh – for surface terms. Compared to the aforementioned works, we here resort to a different hyper-reduction strategy based on empirical quadrature procedure both for volume and surface terms which is  a less intrusive method in its implementation. Furthermore, we rely on a POD-Greedy adaptive sampling strategy based on the definition of an a posteriori error indicator that supports kinematic conditions.

\subsection{Layout of the paper}

The outline of the paper is as follows. In section \ref{sec:formulation}, we present the mathematical formulation of the class of mechanical problems considered in this work (cf. Eq.\eqref{sec:intro:subsec:context:mechanicalPbVar}). In section \ref{sec:methodology}, we display our methodology for building the ROM : we first address the solution reproduction problem, and then we extend our approach to the parametric case. Then in section \ref{sec:model}, we present the physical model problem and assess the methodology validity. In section \ref{sec:numerical}, we present numerical investigations for the model problem and, in section \ref{sec:conclusion}, we draw conclusions and outline subjects of ongoing research.

\section{Formulation}\label{sec:formulation}

\subsection{Formulation of the nonlinear quasistatic problem}

We focus on nonlinear small-displacement small-strain mechanical problems with internal variables. We consider the spatial variable $x$ in the Lipschitz domain $\Omega \subset \mathbb{R}^d$ ($d=2$ or $3$), and the time variable $t\in [0, t_{\rm f}]$. We introduce a vector of parameters $\mu$ which belongs to the compact $\mathcal{P}\subset \mathbb{R}^P$, where $P$ is the number of parameters. As already mentioned, the vector $\mu$ can contain physical parameters (coefficients of the constitutive equations), or geometrical parameters of the problem. We denote by $u$ the primal variable of the mechanical problem (displacement), and we denote by $\mathcal{X}$ the Hilbert space to which the field $u$ belongs. The constitutive equations are assumed to be nonlinear. The system depends implicitly on the displacements history by one implicit differential equation (which includes nonlinear behaviours such as elastoplasticity or viscoplasticity). In this framework, the description of the mechanical state boils down to the knowledge of the displacement field ($u_\mu$), the stress field (the Cauchy tensor $\sigma_\mu$) and the internal variables ($\gamma_\mu$). In this work, we address only quasi-static formulations, which means that we omit the inertial term from the equilibrium equations. Time evolution is described by the system of ordinary differential equations in each point in $\Omega$:

\begin{equation}
\label{sec:intro:subsec:context:mechanicalPb:continuous}
\left\{
\begin{array}{rcl}
- \nabla \cdot \sigma_\mu &=& f_v\\
\sigma_\mu &=& \mathcal{F}^{\sigma}_\mu\left(\nabla_s u_\mu, \ \gamma_\mu\right) \\
\dot{\gamma}_\mu &=& \mathcal{F}^{\gamma}_\mu\left(\sigma_\mu, \gamma_\mu\right)
\end{array}
\right.\quad \text{ + Boundary Conditions (BCs)}
\end{equation}

\noindent where the nonlinear operator $\mathcal{F}^{\sigma}_\mu$ stands for the constitutive equation that maps the state of stresses in the material from the knowledge of deformations ($\nabla_s$ is the symmetric part of the gradient, $\nabla_s\cdot = \frac{1}{2}(\nabla \cdot  + \nabla^T\cdot )$) and internal variables, while the nonlinear operator $\mathcal{F}^{\gamma}_\mu$ denotes an equation of evolution of internal variables within the material. The first equation in the system below describes the equilibrium of our system. The boundary conditions that we consider in this contribution will be detailed later. In this paper, we consider situations where the material is not initially preloaded. At the initial time, all fields are assumed to be zero.

In this work, we restrict ourselves to one-time steps time integrators, implying that the knowledge of the mechanical state is derived from the state previously computed and 'ignores' any information from earlier states beyond that provided by the internal variables.. We introduce the time grid $0=t^{(0)}\leq ..., \leq t^{(K)} = t_{\rm f}$, and we discretize the problem as stated below:

$$u_\mu^{(k)} = u_\mu^{(k-1)} + \Delta u_\mu^{(k)} \quad \text{ and } \quad t^{(k)} = t^{(k-1)} + \Delta t^{(k)}, \qquad \forall k\in \{1,\ldots, K\}$$

\noindent We use a backward Euler discretization scheme for the evolution equation such that the quasi-static discretization of the system boils down to:

\begin{equation}
\label{sec:intro:subsec:context:mechanicalPb:quasistatic}
\left\{
\begin{array}{rclrcl}
-\nabla \cdot \sigma^{(k)}_\mu &=& f_v^{(k)} \qquad &\text{on}&\quad  \ \Omega \\
\gamma^{(k)}_\mu &=& \gamma^{(k-1)}_\mu + \Delta t^{(k)}\mathcal{F}_{\mu}^{\gamma}\left(\sigma^{(k)}_\mu, \gamma^{(k)}_\mu\right)\qquad &\text{on}&\quad  \ \Omega  \\
\sigma^{(k)}_\mu &=& \mathcal{F}_{\mu}^{\sigma}\left(\nabla_s u^{(k)}_\mu, \ \gamma^{(k)}_\mu \right)\qquad &\text{on}&\quad  \ \Omega  \\
&+& \text{BCs} & & \\
\end{array}
\right.
\end{equation}

\noindent Theoretically, stresses can be considered as internal variables. We choose for convenience to restate the problem by displaying only the stress variable in our formulation, as follows:

\begin{equation}
\label{sec:intro:subsec:context:mechanicalPb}
\left\{
\begin{array}{rclrcl}
-\nabla \cdot \sigma^{(k)}_\mu &=& f_v^{(k)} \qquad &\text{on}&\quad  \ \Omega \\
\sigma^{(k)}_\mu &=& \mathcal{F}_{\mu}^{(k)}\left(u^{(k)}_\mu, \ u^{(k-1)}_\mu,  \sigma^{(k-1)}_\mu \right)\qquad &\text{on}&\quad  \ \Omega  \\
&+& \text{BCs} & & \\
\end{array}
\right.
\end{equation}

\noindent where $\mathcal{F}_\mu\left(.,.\right)$ is an appropriate nonlinear operator. In this framework, internal variables are seen as an inner part of the operator $\mathcal{F}_\mu$. We emphasize that our methodology is appropriate for problems of the form \eqref{sec:intro:subsec:context:mechanicalPb:quasistatic}, although we further define it for problems of the form \eqref{sec:intro:subsec:context:mechanicalPb}.

In our study, we consider both non-homogeneous Neumann conditions and homogeneous Dirichlet conditions for suitable linear combinations of the state variables. We assume that the displacement field belongs to the kernel of this form. This choice enables us to model arbitrary linear relations on the displacement field. Other than homogeneous Dirichlet conditions, it supports for instance uniform translation of unknown amplitude of a subpart, or any other arbitrary linear relation between the displacement degrees of freedoms (DOFs) of $\Omega$ accounting for kinematic links between subparts of the system. Such boundary conditions are expressed as:

\begin{equation}
\label{sec:intro:subsec:context:mechanicalPb:boundaryconditions}
\left\{
\begin{array}{rclrcl}
\sigma^{(k)}_\mu \cdot n &=& f_s^{(k)}\qquad &\text{on}&\quad  \ \Gamma_n  \\
c(u_\mu^{(k)}) &=& 0\qquad &\text{on}&\quad  \ \Omega
\end{array}
\right.
\end{equation}

\noindent where $n$ is the outward normal to the boundary $\Gamma_n$, and $f_v^{(k)}$ (resp.$f_s^{(k)}$ ) is the volumic (resp. surfacic) force applied to the system, and $c$ the previously mentionned linear form.
The variational form of the equilibrium equation given by Eq.\eqref{sec:intro:subsec:context:mechanicalPb:continuous} reduces to the following residual expression:

\begin{equation}
\label{sec:intro:subsec:context:mechanicalPbVar:ExpressResidual}
\mathcal{R}_\mu^\sigma \left(\sigma_\mu^{(k)},\ v\right) = \int_{\Omega} \sigma^{(k)}_\mu:\varepsilon(v) \ dx - \int_{\Omega} f_v v \ dx -  \int_{\Gamma_n} f_s v \ ds, \qquad \forall v\in \mathcal{X} 
\end{equation}

\noindent Finally, the variational problem investigated in this contribution  can be summarized as follows:

\begin{equation}
\label{sec:intro:subsec:context:mechanicalPbVar}
\forall k\in \{1,..., K\}, \ \text{Find}\ u^{(k)}_\mu\in \mathcal{X}_{\rm bc}\ \text{s.t.} \ , \quad 
\left\{
\begin{array}{rclrcl}
&\mathcal{R}_\mu\left(u_\mu^{(k)}, \ \ u_\mu^{(k-1)}, \ \sigma_\mu^{(k-1)}, \ v\right) = 0&, \qquad & & \forall v\in \mathcal{X}_{\rm bc} \\
&\sigma^{(k)}_\mu = \mathcal{F}_{\mu}^{(k)}\left(u^{(k)}_\mu, \ u^{(k-1)}_\mu,  \sigma^{(k-1)}_\mu \right)&\qquad &\text{on}&\quad  \ \Omega  \\
&\sigma^{(k)}_\mu \cdot n = f_s^{(k)}&\qquad &\text{on}&\quad  \ \Gamma_n
\end{array}
\right.
\end{equation}

\noindent where $\mathcal{X}_{\rm bc} \coloneq \left\{v\in \mathcal{X}: \ c\left(v\right) = 0, \ \text{on} \ \Omega \right\}$. We also denote :

\begin{equation}
\label{sec:intro:subsec:context:relationexpressionresiduals}
\mathcal{R}_\mu\left(u_\mu^{(k)}, \ u_\mu^{(k-1)}, \ \sigma_\mu^{(k-1)}, \ v\right) = \mathcal{R}^{\sigma}_\mu \left(\mathcal{F}_{\mu}^{(k)}\left(u^{(k)}_\mu, \ u^{(k-1)}_\mu,  \sigma^{(k-1)}_\mu \right), \ v\right).
\end{equation}

\subsection{Finite element discretization}

\subsubsection{Notation}

Given the domain $\Omega \subset \mathbb{R}^d$, we consider a HF mesh $\mathcal{T}^{\rm hf} = \left\{\texttt{D}_i\right\}_{i=1}^{N_{\rm e}}$ where $\texttt{D}_1,\ldots,\texttt{D}_{N_{\rm e}}$ are the elements of the mesh, and $N_{\rm e}$ denotes the number of elements in the mesh. The $\rm hf$ subscript or superscript stands for HF discretization. We allow ourselves to consider meshes with different types of elements in the same mesh. In particular, in the case studied in this work, we consider meshes that include both three-dimensional (volumic) elements and two-dimensional (surfacic) elements; the integer $N_e$ denotes the total number of volumic and surfacic elements. We refer to meshes with volumic and surfacic element as two-level meshes. Besides, we introduce the continuous Lagrangian finite element (FE) basis $\{\varphi_i\}_{i=1}^{\mathcal{N}_{\rm no}}$ associated with the mesh $\mathcal{T}^{\rm hf}$, whose number of nodes is $\mathcal{N}_{\rm no}$. The FE space for the primal unknown is thus defined as follows:

\begin{equation}
\mathcal{X}^{\rm hf}\coloneqq \text{span}\left\{\varphi_i e_j, \ i\in\{1,..., \mathcal{N}_{\rm no}\}, j\in \{1,...,d\} \right\}
\end{equation}

\noindent where $e_1, ..., e_d$ are the vectors of the canonical basis. We further define the nodes $\{x^{\rm hf, \rm no}_i\}_{i=1}^{\mathcal{N}_{\rm no}}$, the quadrature points $\{x^{\rm hf, \rm qd}_i\}_{i=1}^{\mathcal{N}_{\rm qd}}$ associated to the HF-mesh and to the FE discretization and the application  $T^{\rm hf, \rm no}$ (resp. $T^{\rm hf, \rm qd}$ for the quadrature points) which links the global indexing of the DOFs (resp. unknowns at quadrature points) of the HF-mesh to the local indexing of a specified element. The $i_{\rm loc}$-th DOF in the $q$-th element local indexing is associated to the $i_{\rm glob}$ DOF in the global indexing:
 
 $$T^{\rm hf,\rm no}\left(i_{\rm loc}, q\right) = i_{\rm glob}, \quad  i_{\rm loc}=1,\ldots, n_{\rm lp}^q \ \text{and} \ q=1,\ldots,N_{\rm e}$$

\noindent where $n_{\rm lp}^q$ is the number of DOFs in the $q$-th element of the mesh. To further clarify the notation, we denote by $\mathbf{u}\in \mathbb{R}^{\mathcal{N}}$ the FE discrete vector of displacements and $\bm{\sigma}\in \mathbb{R}^{\mathcal{N}_g}$ the stress counterpart, where $\mathcal{N}=d\mathcal{N}_{\rm no}$ is the dimension of the space $\mathcal{X}^{\rm hf}$ and $\mathcal{N}_g = \frac{d(d+1)}{2}\mathcal{N}_{\rm qd}$.

If the mesh contains a single type of element, $T^{\rm hf,\rm no}$ is the connectivity matrix. In the perspective of a hyper-reduced formulation, we introduce two elemental restriction operators: the nodal elemental restriction operators $\{ \mathbf{E}_q^{\rm no} \}_{q=1}^{N_{\rm e}}$ which restrict fields defined at nodes to the q-th element of the mesh (e.g. for displacements) and the quadrature restriction operators $\{ \mathbf{E}_q^{\rm qd} \}_{q=1}^{N_{\rm e}}$ which deals with fields defined at quadrature nodes (e.g. for stresses and internal variables):

$$\left(\mathbf{E}^{\rm no}_q \mathbf{u} \right)_{i_{\rm loc}} = \left(\mathbf{u}\right)_{T^{\rm hf,\rm no}\left(i_{\rm loc}, q\right)}\quad \text{ and } \quad \left(\mathbf{E}^{\rm qd}_{q'} \bm{\sigma} \right)_{j_{\rm loc}} =  \left(\bm{\sigma}\right)_{T^{\rm hf,\rm qd}\left(j_{\rm loc}, q'\right)} $$

\noindent If the restriction refers to a subpart of the mesh, a superscript on the restriction operator is added (for instance, $\mathbf{E}_q^{{\rm no}, \Gamma_n}$ for the nodes of the boundary elements).

\subsubsection{Formulation}

We denote by $\{\mathbf{u}_\mu^{{\rm hf}, (k)}\}_{k=1}^K$ the FE approximation of the displacement (primal variable) given by the HF-model at all times, whereas $\{\bm{\sigma}_\mu^{{\rm hf}, (k)}\}_{k=1}^K$ stand for the stress tensor fields. We state the finite element discretization of Eq.\eqref{sec:intro:subsec:context:mechanicalPbVar}:

\begin{equation}
\label{sec:form:subsec:fe:eq:discritsystem}
\forall k\in \{1,..., K\}, \ \text{Find}\ \mathbf{u}^{{\rm hf}, (k)}_\mu\in \mathcal{X}_{\rm bc}^{\rm hf}\ \text{s.t.} \ , \quad \left\{
\begin{array}{rcl}
&\mathcal{R}^{\rm hf}_\mu\left(\mathbf{u}^{{\rm hf}, (k)}_\mu, \ \mathbf{u}^{{\rm hf}, (k-1)}_\mu, \ \bm{\sigma}^{{\rm hf}, (k-1)}_\mu, \mathbf{v}\right)=0& \qquad \forall \mathbf{v} \in \mathcal{X}_{\rm bc}^{\rm hf}\\
&\bm{\sigma}^{{\rm hf}, (k)}_\mu  = \mathcal{F}^{\rm hf}_\mu\left(\mathbf{u}^{{\rm hf}, (k)}_\mu, \mathbf{u}^{{\rm hf}, (k-1)}_\mu, \ \bm{\sigma}^{{\rm hf}, (k-1)}_\mu\right)&\\
\end{array}
\right.
\end{equation}

\noindent where $\mathcal{X}_{\rm bc}^{\rm hf} \coloneq \left\{\mathbf{v}\in \mathcal{X}^{\rm hf}: \quad  \mathbf{B}\mathbf{v} = 0\right\}$ depicts the test space for displacements, and $\mathbf{B}\in \mathbb{R}^{\mathcal{N}_d\times\mathcal{N}}$ is the kinematic relationship matrix. $\mathcal{N}_d$ stands for the number of linear relations between degrees of freedom that we intend to enforce. Such a formulation on the boundary conditions implies that the kinematic linear application depends neither on time nor on the parameter. The operators $\mathcal{R}^{\rm hf}_\mu$ and $\mathcal{F}^{\rm hf}_\mu$ stands for the discrete counterparts of the continuous operators $\mathcal{R}_\mu$ and $\mathcal{F}_\mu$ introduced in Eq.\eqref{sec:intro:subsec:context:mechanicalPbVar}. Besides, we recall that the displacements are unknowns determined at the nodes of the mesh whereas the discrete stresses are vectors defined at the quadrature points. In practice, the finite element compute the HF-residuals as sums of elementary contributions:

\begin{equation}
\label{sec:form:subsec:fe:eq:discritsystem:elementwise}
\forall \mathbf{v}\in \mathcal{X}^{\rm hf}, \quad \mathcal{R}^{\rm hf}_\mu\left(\mathbf{u}^{(k)}_\mu, \ \mathbf{u}^{(k-1)}_\mu, \ \bm{\sigma}^{(k-1)}_\mu, \ \mathbf{v}\right)=\sum\limits_{q=1}^{N_e}\mathcal{R}^{\rm hf}_{\mu, q}\left(\mathbf{E}_q^{\rm no}\mathbf{u}^{(k)}_\mu, \ \mathbf{E}_q^{\rm no}\mathbf{u}^{(k-1)}_\mu, \ \mathbf{E}_q^{\rm qd}\bm{\sigma}^{(k-1)}_\mu , \ \mathbf{E}_q^{\rm no}\mathbf{v}\right)
\end{equation}

Since a surfacic force is applied on the boundary (see Eq.\eqref{sec:intro:subsec:context:mechanicalPb:boundaryconditions}), two geometric dimensions are involved in our model problem: the domain $\Omega$ and the boundary $\Gamma = \partial \Omega$.  As discussed in the introduction, we here deal with the previously mentioned scenario of a two-level mesh. Within this framework, given \eqref{sec:intro:subsec:context:mechanicalPbVar:ExpressResidual}, the sum expressed in Eq.\eqref{sec:form:subsec:fe:eq:discritsystem:elementwise} is divided into two contributions, one for each geometric dimension:

\begin{align*}
\forall \mathbf{v}\in \mathcal{X}^{\rm hf}, \quad \mathcal{R}^{\rm hf}_\mu\left(\mathbf{u}^{(k)}_\mu, \ \mathbf{u}^{(k-1)}_\mu, \ \bm{\sigma}^{(k-1)}_\mu, \ \mathbf{v}\right)&=\sum\limits_{q=1}^{N_e^{\Omega}}\mathcal{R}^{\rm hf}_{\mu, \Omega_q}\left(\mathbf{E}_q^{{\rm no}, \Omega}\mathbf{u}^{(k)}_\mu, \ \mathbf{E}_q^{{\rm no}, \Omega}\mathbf{u}^{(k-1)}_\mu, \ \mathbf{E}_q^{{\rm qd}, \Omega}\bm{\sigma}^{(k-1)}_\mu, \ \mathbf{E}_q^{{\rm no}, \Omega}\mathbf{v} \right) \\
&+ \sum\limits_{q'=1}^{N_e^{\Gamma_n}}\mathcal{R}^{\rm hf}_{\mu, \Gamma_{n,q'}}\left(\mathbf{E}_{q'}^{{\rm no}, \Gamma_n}\mathbf{u}^{(k)}_\mu, \ \mathbf{E}_{q'}^{{\rm no}, \Gamma_n}\mathbf{u}^{(k-1)}_\mu, \ \mathbf{E}_{q'}^{{\rm qd}, \Gamma_n}\bm{\sigma}^{(k-1)}_\mu, \ \mathbf{E}_{q'}^{{\rm no}, \Gamma_n}\mathbf{v} \right)
\end{align*}

\noindent where we distinguish the residual operators for the integrals over the $q$-th element the domain $\Omega$ ($\mathcal{R}^{\rm hf}_{\Omega_q}$) and the $q'$-th element of the boundary domain $\Gamma_n$ ($\mathcal{R}^{\rm hf}_{\Gamma_{n, q'}}$).

\subsubsection{Dualization of boundary conditions}

So as to comply with the theoretical framework required by the formulations used in our HF setting, the Dirichlet boundary conditions are treated by a dualization of the boundary conditions, namely by introducing Lagrange multipliers. In this setting, the vector solution of the problem at the $k$-th timestep consists of the displacements and the associated Lagrange multipliers $(\mathbf{u}^{(k)}_\mu, \bm{\lambda}^{(k)}_\mu)\in \mathbb{R}^{\mathcal{N}}\times \mathbb{R}^{\mathcal{N}_d}$. The finite element problem assembly amount to the discrete nonlinear system:

\begin{equation}
\label{sec:form:subsec:fe:eq:discritsystem::duallagr}
\forall k \in \{1, \ldots, K \}, \qquad \left\{
\begin{array}{rcl}
\mathbf{R}^{\rm hf}_\mu\left(\mathbf{u}_\mu^{(k)}, \ \mathbf{u}_\mu^{(k-1)}, \ \bm{\sigma}_\mu^{(k-1)}\right) + \mathbf{B}^T \cdot \bm{\lambda}_\mu^{(k)}  &=&0 \\
\mathbf{B}\cdot \mathbf{u}_\mu^{(k)} &=& 0
\end{array}
\right.
\end{equation}

We resort to the Newton-Raphson algorithm to solve \eqref{sec:form:subsec:fe:eq:discritsystem::duallagr}. Note that the Jacobian of \eqref{sec:form:subsec:fe:eq:discritsystem::duallagr} reads as a saddle point system. Dualization of Dirichlet boundary conditions provides a natural framework to enforce Dirichlet-type conditions in the interior of the domain and/or in points that do not coincide with the nodes of the mesh. We provide further details on the solution algorithm in Appendix \ref{app:newton}.

\section{Methodology}\label{sec:methodology}

We seek the reduced-order solution as a linear combination of modes:

\begin{equation}
\widehat{\mathbf{u}}^{(k)}_\mu = \sum\limits_{n=1}^{N_u} \left(\widehat{\mathbf{\alpha}}^{(k)}_{u,\mu}\right)_n\bm{\zeta}_{u,n} = \mathbf{Z}_u\widehat{\mathbf{\alpha}}^{(k)}_{u,\mu}
\end{equation}

\noindent where $\widehat{\mathbf{\alpha}}^{(k)}_{u,\mu}\in \mathbb{R}^{N_u}$ are referred to as generalized coordinates and $\mathcal{Z}_{N_u} = \text{span}\left\{\zeta_{u,n} \right\}$ is the primal reduced space.  The Galerkin ROM is obtained by projecting the discrete residual operator (onto the Eq.\eqref{sec:form:subsec:fe:eq:discritsystem}) onto the primal reduced basis. We first consider the situation without Lagrange multipliers for the boundary conditions:

\begin{equation}
\label{eq:ROBSolving}
\mathbf{Z}_u^T\mathbf{R}^{\rm hf}_\mu\left(\widehat{\mathbf{u}}_\mu^{(k)}, \ \widehat{\mathbf{u}}_\mu^{(k-1)}, \ \widehat{\bm{\sigma}}_\mu^{(k-1)}\right)=0
\end{equation}

Since the operator is nonlinear, successive assemblies are required at each iteration of Newton, leading to a bottleneck in terms of computational costs. As stated in section \ref{sec:intro:subsec:purpose}, we develop an hyper-reduction strategy based on an element-wise empirical quadrature in order to tackle this issue. Thus, the knowledge of the empirical quadrature provides a reduced mesh $\mathcal{T}^{\rm red}$. Assembling the ROM on this mesh speeds up CPU time for the online step. Towards this end, we define the indices associated with the 'sampled' elements. For example, for a two-level discretizations (formulation in Eq.\eqref{sec:form:subsec:fediscrit:eq:2levelmesh}), we have two subsets $I_{\rm eq}^{\Omega} \subset \{1, ..., N^\Omega_e\}$ and $I_{\rm eq}^{\Gamma_n} \subset \{1, ..., N^{\Gamma_n}_e\}$ such that:

\begin{equation}
\label{sec:form:subsec:fediscrit:eq:2levelmesh}
\begin{array}{rcl}
\forall \mathbf{v}\in \mathcal{X}^{\rm hf}, \quad &&\mathcal{R}_{\mu}^{\rm eq}\left(\mathbf{u}^{(k)}, \ \mathbf{u}^{(k-1)}, \ \bm{\sigma}^{(k-1)}, \mathbf{v}\right) \\ \\ &=& \mathcal{R}_{\mu, \Omega}^{\rm eq}\left(\mathbf{u}^{(k)}, \ \mathbf{u}^{(k-1)}, \ \bm{\sigma}^{(k-1)}, \ \mathbf{v}\right) + \mathcal{R}_{\mu, \Gamma}^{\rm eq}\left(\mathbf{u}^{(k)}, \ \mathbf{u}^{(k-1)}, \ \bm{\sigma}^{(k-1)}, \ \mathbf{v}\right) \\ \\
&=&\sum\limits_{q\in I_{\rm eq}^{\Omega} }\rho^{\rm eq, \Omega}_{q}\mathcal{R}^{\rm hf}_{\mu, \Omega_q}\left(\mathbf{E}_q^{{\rm no}, \Omega} \mathbf{u}^{(k)}, \ \mathbf{E}_q^{{\rm no}, \Omega}\mathbf{u}^{(k-1)} , \  \mathbf{E}_q^{{\rm qd}, \Omega}\bm{\sigma}^{(k-1)}, \ \mathbf{E}_q^{{\rm no}, \Omega}\mathbf{v} \right) \\
&+& \sum\limits_{q'\in I_{\rm eq}^{\Gamma_n}}\rho^{\rm eq, \Gamma_n}_{q'}\mathcal{R}^{\rm hf}_{\mu, \Gamma_{n, q'}}\left(\mathbf{E}_{q'}^{{\rm no}, \Gamma_n} \mathbf{u}^{(k)}, \ \mathbf{E}_{q'}^{{\rm no}, \Gamma_n}\mathbf{u}^{(k-1)} , \ \mathbf{E}_{q'}^{{\rm qd}, \Gamma_n}\bm{\sigma}^{(k-1)}, \ \mathbf{E}_{q'}^{{\rm no}, \Gamma_n}\mathbf{v} \right)
\end{array}
\end{equation}

\noindent where $\bm{\rho}^{\rm eq, *} = [\rho^{\rm eq, *}_1,...,  \rho^{\rm eq, *}_{N^*_e}]\in \mathbb{R}^{N^*_e}$  are sparse vectors of positive weights referred as empirical quadrature rules, where $\rho^{\rm eq, *}_q = 0$ if $q\in I_{\rm eq}^{*}$ for $*\in\{\Omega, \ \Gamma_n\}$. It is sufficient to have access to these sets of indices to produce a reduced mesh by considering only the cells with non-zero weights.

Furthermore, the analysis of the mechanical state of the system implies the knowledge of the stresses within the material. To this end, we decompose similarly the stress on an related reduced space $\mathcal{Z}_{N_\sigma}$. This stress basis will namely be used to define an \emph{a posteriori} error indicator.

\subsection{Solution reproduction problem\label{sec:metho:subsec:srpb}}

At first, we omit the parametric variability. In this section, we provide the strategy for constructing a reduced basis and reduced mesh thanks to an empirical quadrature. Our objective here is to reproduce the result obtained in a HF simulation through our reduced problem. The solution reproduction problem is of limited interest; nonetheless, it remains the necessary initial step towards the implementation of an efficient ROM for the parametric problem. The treatment of this sub-problem allows both the design of blocks of algorithms that can be easily reused in the parametric framework, and the provision of validation tests for the latter. This approach is divided into two steps: an \emph{offline} phase where we build reduced bases (displacement and stress), and a reduced mesh, then an \emph{online} phase, which consists in computing the generalized coefficients for both the displacement and the stress. The computation of the coefficients for the stress involves an additional processing with a Gappy-POD procedure (section \ref{subsec:srpb:subsubsec:gappy}). At this point in our study, that is to say without taking into account the design of an error indicator, our reduced model is made up of two reduced bases, one empirical quadrature rule (if not two  in the case of a two-level mesh) and a reduced mesh.

\begin{algorithm}[h]
\caption{Solution Reproduction Problem}\label{alg:srpb}
\begin{algorithmic}
\State \textbf{Online step}
\State Compute the HF-snapshots \Comment{Call of \textsf{code$\_$aster}}
\State Construction of the reduced order basis ($\mathcal{Z}_{N_u}$ and $\mathcal{Z}_{N_\sigma}$)\Comment{Section \ref{subsec:srpb:subsubsec:datacomp}}
\State Empirical Quadrature procedure $\mathbf{\rho}^{\rm eq}$ \Comment{Section \ref{subsec:srpb:subsubsec:eq}}
\State \textbf{Offline step}
\State Compute the primal generalized coordinates $\left\{\mathbf{\widehat{\alpha}}^{(k)}_{u, \mu}\right\}_{k=1}^K$ \Comment{Equation \eqref{eq:ROBSolving}}
\State Compute the dual generalized coordinates using Gappy-POD $\left\{\mathbf{\widehat{\alpha}}^{(k)}_{\sigma, \mu}\right\}_{k=1}^K$\Comment{Section \ref{subsec:srpb:subsubsec:gappy}}
\end{algorithmic}
\end{algorithm}

\subsubsection{Data compression using Proper Orthogonal Decomposition}\label{subsec:srpb:subsubsec:datacomp}

We resort to the method of snapshots\cite{sirovich1987turbulence} to generate both reduced order bases (ROB). We discuss the methodology for the case of the displacement variable; in the case of the stress variable is treated in a similar way. We define the Gramian matrix $\mathbf{C}\in \mathbb{R}^{K\times K}$ associated to a given scalar product $(\mathbf{C}_u)_{i,j} = (u^{\rm hf, (i)}_\mu, u^{\rm hf, (j)}_\mu) = (\mathbf{u}^{\rm hf, (j)}_\mu)^T\bm{X}_u \mathbf{u}^{\rm hf, (i)}_\mu$. Then, we solve the eigenvalue problem:

\begin{equation}
\mathbf{C}_u\varphi_n = \lambda_n \varphi_n, \quad \lambda_1 \geq ...\geq \lambda_K \geq 0
\end{equation}

\noindent to obtain the eigenpairs $(\lambda_n, \varphi_n)$ for $n=1\dotsb N_u$. The number of selected POD modes is chosen according to the following energy criterion: 

\begin{equation}
\label{subsec:srpb:subsubsec:datacomp:eq:energycrit}
N_u = \min \left\{Q\in \mathbb{N}, \quad \sum\limits_{q=1}^Q \lambda_q \geq \left(1 - \varepsilon_{\rm POD, u}^2\right) \sum\limits_{q=1}^K \lambda_q \right\}
\end{equation}

\noindent where $\varepsilon_{\rm POD, u}$ is a user-defined tolerance. It is then now possible to define the POD modes, which will provide the reduced basis for the displacements:

\begin{equation}
\bm{\zeta}_{u,n} = \frac{1}{\sqrt{\lambda_n}}\sum\limits_{k=1}^K \left(\varphi_n\right)_k \mathbf{u}_k
\end{equation}

\noindent In conclusion, given the snapshots $\{u^{\rm hf, (k)}_\mu \}_{k=1}^K$, a scalar product $(.,.)$ and the tolerance $\varepsilon_{\rm POD, u}$, the POD procedure returns the reduced order basis:

\begin{equation}
\mathbf{Z}_u = \text{POD}\left\{\left\{u^{\rm hf, (k)}_\mu \right\}_{k=1}^K, \ (.,.), \  \varepsilon_{\rm POD, u}\right\}
\end{equation}

We need to decide on two scalar products in order to carry out the method: one for the displacement field and one for the stress field. For the displacements, a consistent choice would be to consider the $H^1$ norm. One of the limitations in using the industrial code is that we cannot easily retrieve such a matrix. To overcome this issue, we opted for a compression in the sense of an energy norm. More specifically, we consider the energy norm associated with a simpler mechanical case, that of linear elasticity. In so doing, the formulation of the mechanical problem Eq.\eqref{sec:intro:subsec:context:mechanicalPb} (if we omit the time dependence) becomes :

\begin{equation}
\label{sec:methodo:subsec:podgreedy:subsubsec:errorind:mechanicalPbElastic}
\left\{
\begin{array}{rcl}
-\nabla \cdot \sigma_\mu &=& f_v \qquad \text{on}\quad  \ \Omega \\
\sigma_\mu \cdot n &=& f_s\qquad \text{on}\quad  \ \Gamma_n  \\
u_\mu &=& 0\qquad \text{on}\quad  \ \Gamma_d \\
\sigma_\mu &=& \frac{E}{1+\nu}\nabla_s u_\mu + \frac{E}{\left(1+\nu\right)\left(1 - 2\nu\right)} \left(\nabla \cdot u_\mu\right)\mathds{1}
\end{array}
\right.
\end{equation}

\noindent where $E$ is the Young's modulus and $\nu$ is the Poisson coefficient. From a variational point of view, this amounts to considering a case where we are seeking a displacement field $u\in \mathcal{X}^{\rm hf}_{\rm bc}$ such that :

\begin{equation}
a_\mu\left(u,v\right) = F(v)\quad \text{ with } \quad \left\{
\begin{array}{rcl}
a_\mu\left(u,v\right) &=& \displaystyle \int_{\Omega} \frac{E}{1+\nu} \nabla_s u : \nabla_s v + \frac{E}{\left(1+\nu\right)\left(1 - 2\nu\right)} \left(\nabla \cdot u\right) \left(\nabla \cdot v\right) \ dx \\ \\
F(v) &=& \displaystyle \int_{\Omega} f_v v + \int_{\Gamma_n} f_s v  \\
\end{array}
\right.
\end{equation}

The $a_\mu : \mathcal{X} \rightarrow \mathbb{R}$ is a symmetric, coercive and continuous bilinear form. As a consequence of Korn and Poincaré's inequalities, it defines an equivalent norm of $H^1$: $\forall w\in \mathcal{X}, \ \norm{w}_{a_\mu} = \sqrt{a_\mu\left(u, v\right)}$. However, this energy norm is parametric. To circumvent this issue, we chose to consider the energy norm for the centroid of the parameters $\overline{\mu}\in \mathcal{P}$: $\mathbf{X}_u = \mathbf{K}_{\overline{\mu}}$, where $ \mathbf{K}_{\overline{\mu}}$ is the stiffness matrix obtained for an elastic problem and the vector of parameters $\overline{\mu}$ (or at the components of the vector corresponding to the elastic behaviour). As for the stress field, we consider as the scalar product matrix the diagonal matrix of the HF quadrature weights: $\mathbf{X}_\sigma = \text{diag}(\rho^{\rm hf}_1, \dotsb, \rho^{\rm hf}_{\mathcal{N}_g})$.

\begin{remark}
We have chosen not to compress simultaneously the displacements and the constraints. This decision is motivated by the different roles of both variables in our problem. The problem formulation only involves the displacement field. Therefore, a global compression would require a reorthonormalization of the displacement modes. There is no guarantee of a bijection between these new modes and the global modes. The use of two independent bases for displacement and stress fields helps circumvent this issue.
\end{remark}

\subsubsection{Hyper-reduction via empirical quadrature procedures}\label{subsec:srpb:subsubsec:eq}

In this section, we aim at finding $\mathcal{R}^{\rm eq}_{\mu, \Omega}$ and $\mathcal{R}^{\rm eq}_{\mu, \Gamma}$ according to the separation of the residual described in Eq.\eqref{sec:form:subsec:fediscrit:eq:2levelmesh}. This is done by two distinct calls to the hyper-reduction process described hereafter, one for each level of the mesh. For a given level, the objective of the procedure is to provide an empirical residual defined from the empirical quadrature rule $\bm{\rho}^{\rm eq}$ as given below:

\begin{equation}
\mathcal{R}^{\sigma, \rm eq}\left(\bm{\sigma}^{(k)}, \ \mathbf{v}\right)=\sum\limits_{q\in I_{\rm eq}}\rho^{\rm eq}_{q}\mathcal{R}^{\sigma, \rm hf}_{q}\left(\mathbf{E}_q^{\rm qd}\bm{\sigma}^{(k)}, \ \mathbf{E}_q^{\rm no}\mathbf{v} \right), \qquad \forall \mathbf{v}\in \mathcal{X}^{\rm hf}_{\rm bc}
\end{equation}

\noindent where $\{\bm{\sigma}^{(k)}\}_{k=1}^K$ are the HF snapshots of the problem and $R^{\sigma}$ is defined in Eq.\eqref{sec:intro:subsec:context:relationexpressionresiduals}. In the online phase, the solution is sought on the primal reduced space $\mathcal{Z}_{N_u}$. Therefore, it is sufficient to have a good approximation of the residual on the space spanned by the reduced order basis vector $\text{span}\{\bm{\zeta}_{u,n}\}\subset \mathcal{X}^{\rm hf}$.

Given a tolerance $\delta > 0$, the empirical quadrature rule $\bm{\rho}^{\rm eq}$ should satisfy the following conditions:
\begin{enumerate}
    \item the number of nonzero entries in $\bm{\rho}^{\rm eq}$ should be as small as possible,
    \item the entries of $\bm{\rho}^{\rm eq}$ should be non-negative,
    \item (\emph{constant-function constraints}) the measure of the domain should be conserved: 
    $$\left|\sum\limits_{q=1}^{N_e} \rho^{\rm eq}_q \left|K_q\right| - \left|\Omega\right|\right| < \delta \left|\Omega\right| $$
    \item (\emph{manifold accuracy constraints}) the empirical and HF residuals should be close, meaning that for every primal mode $\bm{\zeta}_{u,n}$ and HF snapshot $\left(\mathbf{u}^{(k)} , \bm{\sigma}^{(k)}\right)$, we have:
\end{enumerate}

\begin{equation}
\left|\sum\limits_{q\in I_{\rm eq}}\rho^{\rm eq}_{q}\mathcal{R}^{\sigma, \rm hf}_{q}\left(\mathbf{E}_q^{\rm qd}\bm{\sigma}^{(k)}, \ \mathbf{E}_q^{\rm no}\bm{\zeta}_{u,n} \right) - \mathcal{R}^{\sigma, \rm hf}\left(\bm{\sigma}^{(k)}, \ \bm{\zeta}_{u,n} \right)\right| \leq \delta \left| \mathcal{R}^{\sigma, \rm hf}\left(\bm{\sigma}^{(k)}, \ \bm{\zeta}_{u,n} \right)\right|
\end{equation}

All these constraints enable us to recast the empirical quadrature problem as a $\ell_0$ pseudo-norm minimisation problem, known as the sparse representation problem:

\begin{equation}
\min\limits_{\bm{\rho}\in N_e} \norm{\bm{\rho} }_{\ell_0}\quad \text{s.t.}\quad
\left\{\begin{array}{rcl}
\norm{\mathbf{G}\bm{\rho} - \mathbf{y}}_* &\leq& \delta \norm{\mathbf{y}}_* \\
\bm{\rho} & \geq & 0
\end{array}
\right.
\end{equation}
for a suitable choice of $\mathbf{G}$,$\mathbf{y}$, $\delta$ and $\norm{.}_*$.

The problem is an NP-hard optimization problem (as indicated in Reference \citenum{farhat2014dimensional} citing Reference \citenum{amaldi1998approximability}) and is therefore not directly solvable in practice. Nonetheless, several alternative methods have been devised in the literature, which rely on relaxation methods inspired by signal processing in order to approximate the quadrature rule in polynomial time.

For instance, Reference \citenum{yano2019lp} proposed an approximation which relies on the $\ell_1$ relaxation of the problem where $\norm{.}_* = \norm{.}_{\ell_\infty}$. The relaxed problem can thus be reformulated as a linear programming problem, and solved by resorting to appropriate solvers. Non-negative least squares problems comprise another class of approximation for the sparse representation problem:

\begin{equation}
\bm{\rho}^{\rm eq} = \mathop{\text{argmin}}\limits_{\bm{\rho}\in \mathbb{R}_+^{N_e}}\norm{\mathbf{G}\bm{\rho} - \mathbf{y}}_2 
\end{equation}

As mentioned in section \ref{sec1}, hyper-reductions methods founded on non-orthogonal matching pursuit algorithms\cite{mallat1993matching}\cite{yaghoobi2015fast} have been developped to this end. Those approaches rely on numerical methods for sparse inexact non-negative least-squares initially developped in signal processing. Similarly, Reference \citenum{farhat2014dimensional} suggested a methodology called Energy-Conserving Sampling and Weighting method (ECSW) that was built on Lawson and Hanson's algorithm\cite{lawson1995solving}. This procedure is an active-set method for solving a non-negative least-square problem. The algorithm is modified thanks to an additional stopping criterion, which helps to enforce the sparsity of the solution. Indeed, a criterion on the residuals obtained in the course of the optimization iterations enables to stop the iterations prematurely:

\begin{equation}
\norm{\mathbf{G}\bm{\rho} - \mathbf{y}}_2 \leq \delta \norm{\mathbf{y}}_2 
\end{equation}

\noindent In our work, this inexact least-squares method has been implemented by modifying the routine in the Python module\cite{virtanen2020scipy} $\rm scipy.optimize.nnls$.

\begin{remark}
We comment on the practical implementation of our approach in this setting in comparison with previous work on hyper-reduction processes in the scope of mechanical problems with internal variables. Indeed, from the expressions given in Eqs \eqref{sec:intro:subsec:context:mechanicalPb:quasistatic} and \eqref{sec:intro:subsec:context:mechanicalPbVar:ExpressResidual}, we may notice that the residual is an operator that explicitly depends on the mechanical state at the current and previous time step. Therefore, this mechanical state includes the internal variables. In previous works\cite{zahr2017multilevel, farhat2014dimensional, iollo2022adaptive}, the exact residual operators used for finite element calculations were called by the hyperreduction process. Therefore, it is required to have explicit information about the internal variables. Several strategies are then available to address this issue explicitly. Reference \citenum{farhat2014dimensional} proposes a storage of the internal variables in addition to the knowledge of the displacement fields to estimate the problem. Another technique proposed in Reference \citenum{iollo2022adaptive} is to call the reduced solver with an exact quadrature rule to gather estimates of the internal variables. Our strategy is slightly different, owing to the inherent restrictions of applying the methodology to an industrial setting and trying to be as non-intrusive as possible. In order to comply with the underlying technical restriction of the industrial code in use, we propose a slightly different strategy. We have opted to reconstruct the integrals of the variational formulation expressed in Eq.\eqref{sec:intro:subsec:context:mechanicalPbVar} outside the HF code and to use these recalculated elementary integrals for the hyper-reduction operation. This entails keeping the calculations performed outside the fidelity code to a bare minimum so that the data used for the learning process is as close as possible to the calculations performed in the assembly in HF practice. Therefore, we extract from the HF code the stress fields at the integration points, the HF gradients of the displacement modes at the integration points (call to the industrial code) and the HF quadrature rule on all the integration points of the mesh. The knowledge of these fields then enables us to reconstruct the said integrals and to carry out the hyper-reduction processes. We may notice that this approach is founded upon the variational formulation, and thus does not demand the knowledge of the internal variables because the information is contained within the stress field itself. However, this procedure involves a slightly higher memory storage cost and an additional call to the HF code for each calculation of the displacement reduced basis (to derive the gradients of the modes at the quadrature points). Nevertheless, this cost remains negligible with respect to a call to the HF code for a complete calculation. Note that the implementation strategy used here is restricted to problems of the form  \eqref{sec:intro:subsec:context:mechanicalPbVar}.
\end{remark}

%\begin{remark}
%Within the literature on hyper-reduction methods, and in particular the ECSW method and its extensions, work has been based on robust black boxes implemented in matlab toolboxes\cite{jain2018simulation}\cite{iollo2022adaptive}, in particular the \verb!lsqnoneg! routine which implements the algorithm presented in the reference \cite{iollo2022adaptive}. In the framework of our work founded on open-source softwares, and interfaced with Python, we extracted the implementation of Lawson's algorithm, in order to reimplement the modified active set algorithm in a Python environment.
%\end{remark}

\subsubsection{Reconstruction of the stress by Gappy-POD}\label{subsec:srpb:subsubsec:gappy}

At the end of a call to the reduced solver, we obtain the reduced solutions in displacement and the related stress (by integration of the constitutive law) at the sampled elements by the empirical quadrature. Nevertheless, these stress vectors do not belong, without loss of generality, to the reduced space designed for the stresses. Indeed, these constraints are derived from the integration of the constitutive law in the HF code from the knowledge of the reduced solution in displacement. When using reduced meshes, the information about the stress is restricted to the quadrature points of the sampled elements. Yet, the description of the mechanical state requires the knowledge of the stress field on the HF mesh. It is thus essential to reconstruct the field on the entire mesh. Furthermore, even without any hyper-reduction procedure, the stresses obtained have no reason to belong to the earlier produced reduced basis, even though no reduced mesh is used. This arises from generating both reduces bases independently. In order to overcome both challenges, we apply a Gappy-POD algorithm\cite{everson1995karhunen} to determine the generalized coordinates.

\subsubsection{Influence of Lagrange multipliers}\label{sec:metho:subsec:srpb:subsubsec:lagrange}

As previously mentionned, we address arbitrary, homogeneous (right-hand side is 0) kinematic links between the dofs in $\Omega$ in this work, written as $\mathbf{B} \mathbf{u}_\mu^{(k)} = 0$. As will be shown, specifying kinematic links as inputs to the online solver is then no longer necessary, nor is the implementation of any specific treatments for them during online resolution. In this strategy, kinematic links are already taken into account by the reduced basis, which greatly simplifies coding of the online resolution. It is worth noting that such a choice can reduce drastically the number of unknowns and therefore the computational cost. Indeed, we have:

\begin{equation}
\label{eq:ROBSolving_dualBCs}
\begin{array}{rcl}
\mathbf{Z}_u^T\mathbf{R}_{\mu}^{\rm hf}\left(\widehat{\mathbf{u}}^{(k)}_{\mu},  \widehat{\mathbf{u}}^{(k-1)}_{\mu}, \widehat{\bm{\sigma}}^{(k-1)}_{\mu}\right) + \mathbf{Z}_u^T\mathbf{B}^T \bm{\widehat{\lambda}}^{(k)}_{\mu} &=&\mathbf{Z}_u^T\mathbf{R}_{\mu}^{\rm hf}\left(\widehat{\mathbf{u}}^{(k)}_{\mu}, \widehat{\mathbf{u}}^{(k-1)}_{\mu}, \widehat{\bm{\sigma}}^{(k-1)}_{\mu}\right) + \left[\mathbf{B}\mathbf{Z}_u\right]^T \bm{\widehat{\lambda}}^{(k)}_{\mu}\\
&=& \mathbf{Z}_u^T\mathbf{R}_{\mu}^{\rm hf}\left(\widehat{\mathbf{u}}^{(k)}_{\mu}, \widehat{\mathbf{u}}^{(k-1)}_{\mu}, \widehat{\bm{\sigma}}^{(k-1)}_{\mu}\right) =0
\end{array}
\end{equation}

\noindent where we have omitted the stress field so as not to make the equations more cumbersome. By construction, for a given $n$, $\mathbf{\zeta}_{u,n}$ is a linear combination of the snapshots and therefore verify $\mathbf{B}\mathbf{\zeta}_{u,n}=0$. Such a setting reduces the number of unknowns, as the Lagrange multipliers can be ignored. They do not need to appear in the resolution of the nonlinear system, to be stored or to be taken into account in a data compression operation.

\subsection{Adaptive algorithm based on POD-Greedy procedure}

As mentioned previously, we develop an adaptive sampling based on a POD-Greedy strategy. Moreover, we introduce an error indicator correlated to the approximation error, whose evaluation is cost-efficient in terms of computational time. The extension of the above problem to a parametric system raises two challenges: first, the adaptation of the data compression and the empirical quadrature techniques to the iterative process; second, the construction of an error indicator for our model problem. Indeed, we must be able to have some information about the reliability of our  mechanical state estimation for a given parameter.

\begin{algorithm}[h]
\caption{POD-Greedy algorithm}\label{alg:podgreedy}
\begin{algorithmic}[1]
\Require $\Theta_{\rm train}=\{\mu_i\}_{i}^{n_{\rm train}}$, $\varepsilon_{\rm POD, u}$, $\varepsilon_{\rm POD, \sigma}$
\State $\mathcal{Z}_{N_u} = \mathcal{Z}_{N_\sigma} = \emptyset$, $\mu_* = \overline{\mu}$, $\Theta_{*}=\{\mu_*\}$.
\While{ Stop Criterium }
\State Compute $\{\mathbf{u}^{{\rm hf}, (k)}_{\mu^*}\}_{k=1}^K$,  $\{\bm{\sigma}^{{\rm hf}, (k)}_{\mu^*}\}_{k=1}^K$\Comment{Call of \textsf{code$\_$aster}}
\State Compute primal ROB $\mathbf{Z}_u$ \Comment{Section \ref{sec:methodo:subsec:podgreedy:subsubsec:datacompre}}
\State Compute $\bm{\rho}^{\rm eq}$ knowing $\{\zeta_{u,n} \}_{n=1}^{N_u}$ and $\{\bm{\sigma}^{{\rm hf}, (k)}_{\mu}\}_{k\in \{1, .., K\}, \mu\in \Theta_{\rm *}}$ \Comment{Section \ref{subsec:srpb:subsubsec:eq}}
\State Compute the reduced mesh $\mathcal{T}^{\rm red}$
\State Compute dual ROB $\mathbf{Z}_\sigma$ \Comment{Section \ref{sec:methodo:subsec:podgreedy:subsubsec:datacompre}}
\State Construction of $\bm{\Sigma}_N$ for the error indicator
\For{$\mu \in \Theta_{\rm train}$}
\State Solve the ROM for $\mu$ and compute $\Delta_{N, \mu}^{\rm avg}$ \Comment{Section \ref{sec:methodo:subsec:podgreedy:subsubsec:errorind}}
\EndFor
\State $\Delta_N^{\rm av, \max}= \max\limits_{\mu\in \Theta_{\rm train}}\Delta_{N, \mu}^{\rm avg}$
\State $\mu^* = \arg  \max\limits_{\mu\in \Theta_{\rm train}}\Delta_{N, \mu}^{\rm avg}$
\State $\Theta_{\rm *}= \Theta_{\rm *}\cup \{\mu_*\}$
\EndWhile
\end{algorithmic}
\end{algorithm}

\begin{remark}
In the description of the algorithm, we do not state explicitly the stopping criteria for the algorithm. Several options are possible: a tolerance on the minimum value of the residual can be given, or a maximum number of iterations can be imposed by the user. An alternative scenario is to check whether the approximation error for the new parameter to be explored for the current reduced order basis is below a given threshold.  If so, such a criterion illustrates that we have thus already exploited the redundancy of information and the algorithm can stop.
\end{remark}

\begin{figure}[h]
\begin{center}
\includegraphics[scale=1]{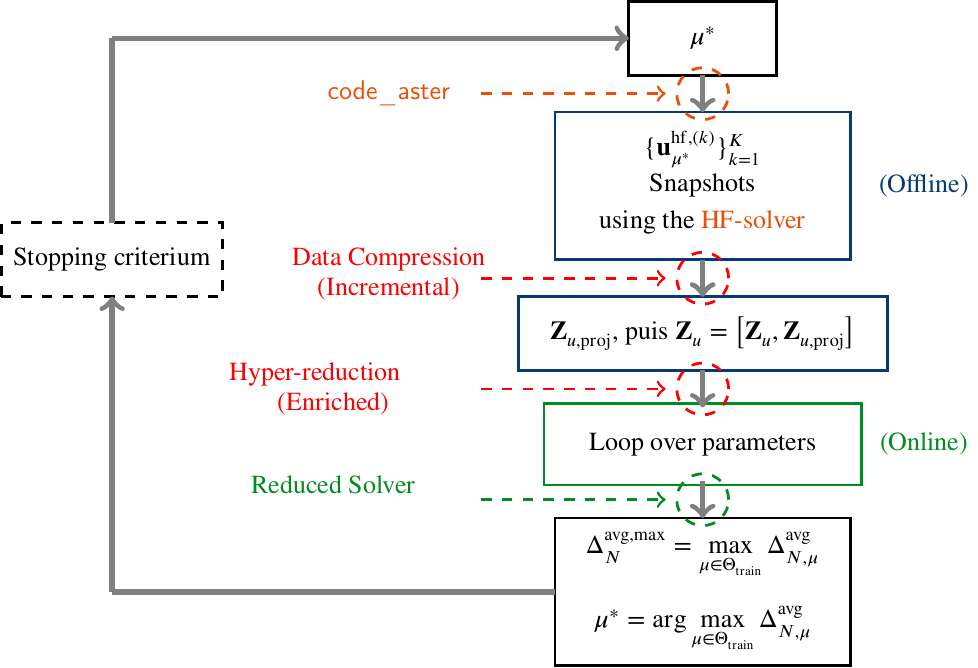} 
\end{center}
\caption{Adaptive algorithm based on POD-Greedy procedure}
\end{figure}

\subsubsection{Data compression}
\label{sec:methodo:subsec:podgreedy:subsubsec:datacompre}

Since we chose to adopt a hierarchical basis, we decided to implement an incremental POD, by applying the POD procedure on the projection of new snapshots on the orthogonal space to the existing basis. Suppose that we have a reduced order basis $\mathbf{Z}_u$ and new HF snapshots $\{u^{\rm hf, (k)}_\mu \}_{k=1}^K$. The new basis is obtained by concatenation:

\begin{equation}
\mathbf{Z}_u = \left[\mathbf{Z}_u, \mathbf{Z}_{\rm proj}\right], \quad \text{where} \quad \mathbf{Z}_{\rm proj}= \text{POD}\left\{\left\{\Pi_{\mathcal{Z}_u^{\bot}, (.,.)} u^{\rm hf, (k)}_\mu \right\}_{k=1}^K, \ (.,.), \  \varepsilon_{\rm POD, u}\right\}
\end{equation}

\noindent where $\Pi_{\mathcal{Z}_u^{\bot}, (.,.)} : \ \mathcal{X}^{\rm hf} \rightarrow \mathcal{Z}_u $ is the orthogonal projection operator onto $\mathcal{Z}_u\subset \mathcal{X}^{\rm hf} $ using the $(.,.)$ scalar product. 

This approach is referred as H-POD in the literature. We refer to Reference \citenum{haasdonk2008reduced} for more details. In terms of memory storage cost, this method does not require to store the eigenvalues between two consecutive iterations. Previous works\cite{iollo2022adaptive} have highlighted the challenge of finding an optimal tolerance for the data compression. Indeed, the compression operator aims to properly describe the snapshots that are supplied to it. Nevertheless, from a numerical standpoint, the vector projection can never be exactly zero. The issue is that the number of selected POD modes is usually chosen based on an energy criterion (cf. Eq.\eqref{subsec:srpb:subsubsec:datacomp:eq:energycrit}). If we apply POD to the projected snapshots, this criterion might be unreliable due to the fact that the energy content of the projected snapshots might be extremely modest if compared to the energy content of the original snapshot set. This observation explains the importance of introducing a criterion based on the relative projection error.  We rely on the regularization approaches given in Reference \citenum{iollo2022adaptive}. The number of modes is chosen according to the following criterion:

\begin{equation}
N^{\rm new} = \min \left\{M \ : \ \max\limits_{k\in\{1,...,K\} }\frac{\norm{\Pi_{\left(\mathcal{Z}_{N_u} \oplus \mathcal{Z}_{u, M}^{\rm new}\right)^{\bot}, (.,.)} u^{\rm hf, (k)}_\mu}}{\norm{u^{\rm hf, (k)}_\mu}} \leq \varepsilon_{\rm POD, u},\quad \mathcal{Z}_{u, M}^{\rm new} = \text{span}\left\{\zeta_{u, m}^{\rm new} \right\}_{m=1}^M
\right\}
\end{equation}

Only the basis vectors that effectively reduce the projection error are added to the reduced basis. The others are treated as noise and are dropped. 
On top of that, for numerical efficiency purposes, we have chosen a criterion prior to the computation of the extra modes. No further POD computation is performed when: :

\begin{equation}
\max\limits_{k\in\{1,...,K\} }\frac{\norm{\Pi_{\left(\mathcal{Z}_{N_u}\right)^{\bot}, (.,.)} u^{\rm hf, (k)}(\mu)}}{\norm{u^{\rm hf, (k)}_\mu}} \leq \varepsilon_{\rm POD, u}
\end{equation}

Based on the very same principle, we assume in this situation that the new snapshots belong to the previously generated reduced space, and there is no update of the basis. This preliminary verification avoids unnecessary offline CPU costs.

\subsubsection{Error indicator}
\label{sec:methodo:subsec:podgreedy:subsubsec:errorind}

We introduce an error indicator to assess the quality of our approximation without having to compute approximation errors, that is to say without having to compute further HF snapshots. We choose to consider a time-averaged error indicator defined as the averaged of the dual norm at each timestep:

\begin{equation}
\Delta_{N, \mu}^{\rm avg} = \sqrt{\frac{1}{K} \sum\limits_{k=1}^K \left( \Delta_{N, \mu}^{(k)} \right)^2}, \quad \text{with} \quad \Delta_{N,\mu}^{(k)} = \sup\limits_{v\in \mathcal{X}^{\rm hf}_{\rm bc}} \frac{\mathcal{R}^{\rm hf}_\mu \left(\widehat{u}^{(k)}_\mu, \widehat{u}^{(k-1)}_\mu ,\widehat{\sigma}^{(k-1)}_\mu, v\right)}{\norm{v}}
\end{equation}

We expect that the error indicator is correlated to the approximation error in displacement prediction (and ideally stress prediction). If so, the error indicator can be used to drive the greedy strategy. We can derive an efficient online/offline strategy which relies on the fact that the stress prediction belongs on a given reduced space $\widehat{\sigma}^{(k)}_\mu  \in \text{span}\{\zeta_{\sigma, n}\}$ and that the reduced residual can be expressed thanks to it.

%\begin{equation}
%\Delta_{N, \mu}^{(k)} = \sup\limits_{v\in \mathcal{X}_{\rm bc}} \frac{\mathcal{R}^{\eq r}\left(\widehat{u}_\mu, v \right)}{\norm{v}} =  \sup\limits_{v\in \mathcal{X}_{\rm bc}} \left[\frac{\int_{\Omega}\widehat{\sigma}_\mu^{(k)}:\varepsilon(v) \ dx - \int_{\Omega} f_v^{(k)}v\ dx  - \int_{\Omega}f_s^{(k)}v \ ds}{\norm{v}}\right]
%\end{equation}

In the following formulation, we assume that the external loadings do not depend on the time variable.  This choice is made for the sake of simplicity, and we can refer to Appendix \ref{app:errorind} for the more general formulation. We introduce the Riesz elements $\psi_n^{\sigma}\in \mathcal{X}^{\rm hf}_{\rm bc}$ associated to the given linear forms:

\begin{equation}
\label{sec:methodo:subsec:podgreedy:subsubsec:errorind:eqdefriesz}
\left(\psi_n^{\sigma}, \ v\right) = \mathcal{L}_n (v), \quad  \forall v\in \mathcal{X}^{\rm hf}_{\rm bc} \quad  \text{ with } \left\{
\begin{array}{rcl}
\mathcal{L}_n(v)&=& \int_{\Omega}\zeta_{\sigma, n}:\varepsilon(v) \ dx, \quad 1 \leq n \leq N_\sigma \\
\mathcal{L}_{N_\sigma + 1} &=& \int_{\Omega} f_v\cdot v \ dx  + \int_{\Gamma_n} f_s \cdot v \ ds 
\end{array}
\right.
\end{equation}

\noindent By means of the decomposition of the stress solution on the reduced basis ($\widehat{\bm{\sigma}}_\mu^{(k)} = \mathbf{Z}_\sigma\widehat{\mathbf{\alpha}}^{(k)}_{\sigma,\mu}$) and the expression of the residual in variational form given by Eq.\refeq{sec:intro:subsec:context:mechanicalPbVar}, we can recast the dual norm calculation as:

\begin{equation}
\Delta_{N, \mu}^{(k)} = \sup\limits_{v\in \mathcal{X}_{\rm bc}} \left[ \sum\limits_{n=1}^{N_\sigma}\left(\mathbf{\widehat{\alpha}}^{(k)}_{\sigma, \mu}\right)_n \frac{\mathcal{L}_n(v)}{\norm{v}} - \frac{\mathcal{L}_{N_\sigma + 1}(v)}{\norm{v}}  \right] = \sup\limits_{v\in \mathcal{X}_{\rm bc}} \frac{\left(\sum\limits_{n=1}^{N_\sigma}\left(\mathbf{\widehat{\alpha}}^{(k)}_{\sigma, \mu}\right)_n \psi_n^\sigma - \psi_{N_\sigma + 1}^\sigma, \ v \right)}{\norm{v}}
\end{equation}

\noindent The dual norm is equal to the norm of its Riesz element, which gives a compact expression for the error indicator:

\begin{equation}
\left(\Delta_{N, \mu}^{(k)} \right)^2 = \norm{\sum\limits_{n=1}^{N_\sigma}\left(\mathbf{\widehat{\alpha}}^{(k)}_{\sigma, \mu}\right)_n \psi_n^\sigma - \psi_{N_\sigma + 1}^\sigma}^2 = \begin{bmatrix}\mathbf{\widehat{\alpha}}^{(k)}_{\sigma, \mu}\\
-1
\end{bmatrix}^T \bm{\Sigma}_N\begin{bmatrix}\mathbf{\widehat{\alpha}}^{(k)}_{\sigma, \mu}\\
-1
\end{bmatrix} = \left(\widetilde{\mathbf{\alpha}}^{(k)}_{\sigma, \mu}\right)^T\bm{\Sigma}_N\widetilde{\mathbf{\alpha}}^{(k)}_{\sigma, \mu}
\end{equation}

\noindent where $\bm{\Sigma}_N \in \mathbb{R}^{N_\sigma + 1, N_\sigma + 1}$ is the Gramian matrix of the Riesz elements previously introduced, i.e $\left(\bm{\Sigma}_N \right)_{n,m}= \left(\psi_n^{\sigma}, \ \psi_m^{\sigma}\right)$, and $\widetilde{\mathbf{\alpha}}^{(k)}_{\sigma, \mu}$ is the concatenation of the generalized coordinates for the stress with $[-1]$.

We shall now discuss the effective calculation of the Riesz elements in the context of the industrial code. In a general manner, these vectors can be determined by solving $N_\sigma + 1$ linear systems defined by Eq.\eqref{sec:methodo:subsec:podgreedy:subsubsec:errorind:eqdefriesz} and will hence fulfil the boundary conditions associated with the system: $\mathbf{B}\bm{\psi}_n^\sigma = 0$. In our framework, it is not straightforward to formulate a problem in variational form by hand or to extract all information to solve Eq.\eqref{sec:methodo:subsec:podgreedy:subsubsec:errorind:eqdefriesz} algebraically outside the FE solver. Nevertheless, functionalities exist to extract Riesz elements of the given linear forms but on a larger space $\mathcal{X}^{\rm hf}$, i.e. for vectors that do not satisfy the boundary conditions of the problem. Indeed, such features are often implemented in industrial codes so that engineers can have access to internal forces vectors or support reaction forces. Such vectors are defined as:

\begin{equation}
\label{sec:methodo:subsec:podgreedy:subsubsec:errorind:discreteRiesz}
\forall \mathbf{v}\in \mathbb{R}^{\mathcal{N}}, \qquad \left\{
\begin{array}{rcl}
\left(\mathbf{F}_n, \ \mathbf{v}\right)_{\ell_2} &=& \int_{\Omega}\zeta_{\sigma, n}:\varepsilon(v) \ dx , \qquad \forall n \in \{1, ..., N_\sigma \} \\ \\
\left(\mathbf{F}_{N_\sigma + 1}, \ \mathbf{v}\right)_{\ell_2} &=& \int_{\Omega} f_v\cdot v \ dx  + \int_{\Gamma_n} f_s \cdot v \ ds 
\end{array}
\right.
\end{equation}

As a reminder, the Riesz elements should belong to the same space as the displacement space. For the sake of consistency, the scalar product used to define them is the scalar product associated to the energy norm for $\overline{\mu}$ (section \ref{subsec:srpb:subsubsec:datacomp}). It is clear from Eq.\eqref{sec:methodo:subsec:podgreedy:subsubsec:errorind:eqdefriesz} that $\mathbf{\psi}_n^\sigma$ is solution to a quadratic optimization problem associated with cost function $\mathbf{v} \mapsto \frac{1}{2} \mathbf{v}^T \mathbf{K}_{\overline{\mu}} \mathbf{v} - \mathbf{v}^T \mathbf{F}_n$ under the equality constraint $\mathbf{B} \mathbf{v} = 0$. The KKT optimality conditions read:

\begin{equation}
\label{sec:methodo:subsec:podgreedy:subsubsec:errorind:KKT}
\left\{
\begin{array}{rcl}
\mathbf{K}_{\overline{\mu}} \bm{\psi}^{\sigma}_n + \mathbf{B}^T\bm{\lambda}&=& \mathbf{F}_n,  \\
\mathbf{B}\bm{\psi}^{\sigma}_n &=& 0
\end{array}
\right.
\end{equation}

\noindent Actually, Eq.\eqref{sec:methodo:subsec:podgreedy:subsubsec:errorind:KKT} defines an easy problem to provide as an input to a FEM solver: it is a linear elastic case for the parameter centroid, with the very same boundary conditions as the HF problem, and an explicit field of nodal forces as an external load (previously computed by Eq.\eqref{sec:methodo:subsec:podgreedy:subsubsec:errorind:discreteRiesz}) It is therefore sufficient to use the HF solver for $N_\sigma + 1$ linear problems. Finally, the parameter-independent matrix that appears in the error indicator definition is computed as follows:

\begin{equation}
\left(\bm{\Sigma}_N\right)_{n,m} = \left(\psi_n^\sigma, \ \psi_m^\sigma \right)_{a_{\overline{\mu}}} = \bm{\psi}_n^\sigma \cdot \left(\mathbf{K}_{\overline{\mu}} \cdot \bm{\psi}_m^\sigma\right) = \left( \bm{\psi}_n^\sigma\right)^T\mathbf{K}_{\overline{\mu}} \cdot \bm{\psi}_m^\sigma, \qquad \forall n,m \in \{ 1, \dotsb, N_{\sigma} + 1 \}
\end{equation}

\begin{remark}
For a given HF FE solver, this strategy is a non-intrusive way to compute the error indicator, since there is no need to retrieve $\mathbf{B}$ or $\mathbf{K}_{\overline{\mu}}$ matrices from the finite element solver.
\end{remark}

\section{Model problem: elastoplastic analysis of a plate with a hole}\label{sec:model}

We validate the approach through the vehicle of an elasto-plastic three-dimensional holed plate. We investigate the problem of a plate with a hole submitted to a traction loading. In this section, we first present the physical formulation of the material constitutive law, then the resolution algorithm used in our work, and finally the details of the configuration used in our numerical example. For the sake of simplicity, we remove the parametric dependence in the notation, which means that the $\mu$ subscript is removed in this section.

\subsection{Elasto-plasticity using a Von Mises Criterion}

\subsubsection{Continuous equations}

We consider a small-displacement small-strain mechanical problem. We assume that the total deformation is the sum of a plastic part ($\varepsilon^{\rm p}$) and an elastic part ($\varepsilon^{\rm el}$):

$$\varepsilon =  \varepsilon^{\rm el} + \varepsilon^{\rm p}$$

\noindent where the plastic deformation comprises the irreversible part of the behavior. The elastic behavior depends on two parameters, the Young's modulus $E$ and the Poisson coefficient $\nu$. The elastic  constitutive equation is :

\begin{equation}
\label{sec:Phys:subsec:Plas:eq:constitutiveElas}
\sigma = \mathcal{F}^{\sigma}\left(\nabla_s u, \varepsilon^p\right) = \frac{E\nu}{\left(1+\nu\right) \left(1-2\nu\right)}\text{Tr}\left(\nabla_s u - \varepsilon^p\right)\mathds{1} + \frac{E}{1+\nu}\left(e - e^p\right)
\end{equation}

\noindent where the deviator of the strain and stress tensors are introduced:

\begin{equation}
\label{sec:Phys:subsec:Plas:eq:constitutiveElas:notations}
e = \text{dev}\left(\nabla_s u\right), \quad e^p = \text{dev}(\varepsilon^p), \quad s = \text{dev}(\sigma), \qquad \text{where} \ \qquad\text{dev}(\tau) = \tau - \frac{1}{3}\text{Tr}(\tau)\mathds{1} 
\end{equation}

We consider a Von Mises criterion for an isotropic hardening. In our analysis, the internal variables that appear in the model are the plastic strain ($\varepsilon^p$) and the cumulative plastic strain ($p$). In the framework of the formulations presented in the previous section, this decision implies that the evolution equations on the internal variables are expressed using the following system:

\begin{equation}
\left(\dot{\varepsilon}^p, \dot{p}\right)= \mathcal{F}^{\gamma}\left(\sigma, \varepsilon^p, p\right) \quad \Leftrightarrow
\left\{
\begin{array}{rcl}
&\sigma^{\rm eq} =  \sqrt{\frac{3}{2} s : s }& \\
&\sigma^{\rm eq } - R(p) \leq  0&\qquad \text{[Von Mises criterion]}\\
&p(t) = \sqrt{\frac{3}{2}}\int_0^t \norm{\dot{\varepsilon}^p\left(\tau\right)} \ d\tau& \\
&\dot{\varepsilon}^p = \dot{p} \frac{3}{2\sigma^{\rm eq }}s \quad  \quad \dot{p}\geq 0  \quad \dot{p}\left[\sigma^{\rm eq} - R(p) \right] = 0 &\qquad  \text{[Normality rule]}
\end{array}
\right.
\end{equation}

\noindent where $\sigma^{\rm eq}$ is an Von Mises equivalent stress and $R(p)$ denotes the elastic limit, and evolves as a function of the cumulative plastic strain $p$. For more insight into the time-discretized formulation, and in the spirit of reproducibility of this research, the reader may find all the details of the numerical procedure in Appendix \ref{app:elastoplasticity}.

\subsection{Physical problem and algorithm}
\label{sec:model:subsec:physpb}

The work hardening curve is chosen to follow a power law (referred as VMIS$\_$ISOT$\_$PUIS in the \textsf{code$\_$aster} database), which implies that the elastic limit evolves on the accumulated plastic strain as follows:

$$R(p) = \sigma_y + \sigma_y\left(\frac{E}{a_{\rm pui}\sigma_y}p\right)^{\frac{1}{n_{\rm pui}}} $$

\noindent where $n$, $a_{\rm pui}$ are strain hardening coefficients and $\sigma_y$ is the initial elastic limit. This algorithm provides us with stable responses on a range of parameters. The resolution procedure used in this work is a elastic predictor-return mapping (plastic corrector)\cite{wilkins1963calculation}. In case of a plastic evolution, the nonlinear equation that ensures the fulfillment of the criterion is solved using the secant method. All the physical parameters of the problem are summarized in Table \ref{sec:model:subsec:physpb:table:physparameters}.

\begin{table}[h!]
\begin{center}
\begin{tabular}{|c|c|c|c|c|}
\hline
$E$ & $\nu$ & $\sigma_y$ & $n_{\rm pui}$ & $a_{\rm pui} $ \\ \hline
MPa & no dim. & MPa & no dim. & no dim. \\ \hline
\end{tabular}
\end{center}
\caption{Summary of the physical parameters\label{sec:model:subsec:physpb:table:physparameters}}
\end{table}

\subsection{Geometric configuration and physical parameters}

We shall study the problem of a three-dimensional plate with a circular hole in its centre and subjected to a tension force. Such a typical example is widely studied in the mechanics literature and is therefore a classical test case for the investigation of algorithms in nonlinear mechanics, namely in elastoplasticity.

\begin{figure}[h]
\begin{center}
\includegraphics[scale=1]{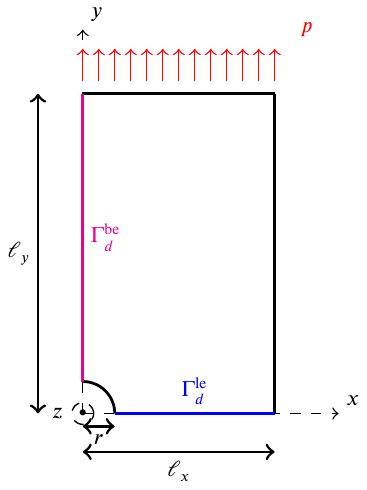} 
\end{center}
\caption{Geometric configuration and loading\label{geom:figure:geometryNumCase}}
\end{figure}

The geometrical domain is narrowed for reasons of symmetry (geometry given on Figure \eqref{geom:figure:geometryNumCase}). We consider that the tension force is only applied on the upper boundary of the plate. We assume that the vertical displacement is homogeneous on the upper boundary, where symmetric boundaries are also applied:

\begin{equation}
\left\{
\begin{array}{rcl}
-\nabla \cdot \sigma = f_v \quad \text{on } \Omega\\
\sigma \cdot e_y = -p \quad \text{on } \Gamma_n^{\rm up}
\end{array}
\right., \qquad \text{ such that } \qquad \left\{
\begin{array}{rcl}
u_y = 0 & \text{on} & \Gamma_d^{\rm be}\\
u_x = 0 \quad &\text{on}& \quad  \Gamma_d^{\rm le}\\
u_z = 0 \quad &\text{on}& \quad  \Gamma_d^{\rm ba} \\
u_z(x_1) &=& u_z(x_2), \quad \forall (x_1, x_2)\in \Gamma_n^{\rm up} \times \Gamma_n^{\rm up} 
\end{array}
\right.
\end{equation}

\noindent where the associated boundaries are defined as:

$$\left\{
\begin{array}{rcl}
\Gamma_d^{\rm be} &=& \left\{y=0, \ \forall (x,z)\in\left[r, \ell_x\right]\times \left[0, \ell_z\right] \right\}\\
 \Gamma_d^{\rm le}& = &\left\{y=0, \ \forall (y,z)\in \left[r, \ell_y\right]\times \left[0, \ell_z\right]\right\}\\
 \Gamma_d^{\rm ba}& =& \left\{y=0, \ \forall (x,y)\in \left[r, \ell_x\right]\times \left[0, \ell_z\right] \setminus \{(x,y), \ \text{ s.t.} \ x^2 + y^2 < r\}\right\}\\
 \Gamma_n^{\rm up} &=&\left\{(x, \ell_y, z), \ \forall x \left[0, \ell_x\right], \forall z\in \left[0, \ell_z\right] \right\}
\end{array}
\right.$$

\section{Numerical results}\label{sec:numerical}

We measure the performance through the previsouly defined energy norm on the FE vectors. We introduce the time-averaged projections errors and approximation errors on the displacement for any $\mu\in \mathcal{P}$:

\begin{equation}
\label{sec:numerical:eq:projerrors}
E_{u, \mu}^{\rm proj, avg} = \frac{\sqrt{\sum\limits_{k=1}^K\norm{\Pi_{\mathcal{Z}_{N_u}^\bot}u^{{\rm hf}, (k)}}^2}}{\sqrt{\sum\limits_{k=1}^K\norm{u^{{\rm hf}, (k)}}^2}} \quad \text{and} \quad E_{u, \mu}^{\rm app, avg} = \frac{\sqrt{\sum\limits_{k=1}^K\norm{u^{{\rm hf}, (k)} - \widehat{u}^{(k)}}^2}}{\sqrt{\sum\limits_{k=1}^K\norm{u^{{\rm hf}, (k)}}^2}}
\end{equation}

For the numerical tests, we treat a strain hardening parameter $a_{\rm pui}$ and the Poisson's ratio $\nu$ as varying parameters (see Table \ref{sec:model:subsec:physpb:table:physparameters}). We define the parameter compact as a Cartesian product of parameter intervals $\mathcal{P} = \mathcal{P}_\nu\times \mathcal{P}_{a_{\rm pui}} = \left[0.21, 0.3\right]\times \left[0.1, 1000\right]$. At last, we introduce the discrete version of this compact $\Theta_{\rm train} = \Theta_{{\rm train}, \nu}\times  \Theta_{{\rm train}, a_{\rm pui}}$, which thus constitutes the training set we shall examine. In order to assess the method, we have carried out numerical tests in several steps. Each step allows to validate specific features of the methodology we have designed. As described in the methodology section, the first step of the validation is the processing of a solution reproduction problem (section \ref{sec:metho:subsec:srpb}), which illustrates the interest of data compression and the construction of a reduced mesh in terms of CPU cost, while maintaining a quality in the approximation of the solution. Afterwards, we shall discuss two parametric cases: afirst, we consider a case with a scalar parameter; second, we consider the case of a two-dimensional parameter.

\begin{figure}[h!]
\begin{center}
\begin{subfigure}[c]{0.3\textwidth}
\begin{center}
\begin{tabular}{|c|c|c|}
\hline
$\left|\Theta_{{\rm train}, \nu}\right|$ &    $\left|\Theta_{{\rm train}, a_{\rm pui}}\right|$ & $K$  \\ \hline
1 & 1 & 20 \\ \hline
\end{tabular}
\end{center}
\caption{\normalsize Solution Reproduction Problem \\ (section \ref{sec:numres:subsec:srpb})}
\end{subfigure}
\quad
\begin{subfigure}[c]{0.3\textwidth}
\begin{center}
\begin{tabular}{|c|c|c|}
\hline
$\left|\Theta_{{\rm train}, \nu}\right|$ &    $\left|\Theta_{{\rm train}, a_{\rm pui}}\right|$ & $K$  \\ \hline
20 & 1 & 10 \\ \hline
\end{tabular}
\end{center}
\caption{Parametric Problem $\mu = \nu$ \\ (section \ref{sec:numres:subsec:param:subsubsec:nu})}
\end{subfigure}
\quad 
\begin{subfigure}[c]{0.3\textwidth}
\begin{center}
\begin{tabular}{|c|c|c|}
\hline
$\left|\Theta_{{\rm train}, \nu}\right|$ &    $\left|\Theta_{{\rm train}, a_{\rm pui}}\right|$ & $K$  \\ \hline
20 & 20 & 10 \\ \hline
\end{tabular}
\end{center}
\caption{Parametric Problem $\mu = \left(\nu, a_{\rm pui}\right)$ \\ (section \ref{sec:numres:subsec:param:subsubsec:nuapui})}
\end{subfigure}
\end{center}
\caption{Summary of the size of the training sets and the number of timesteps ($K$) used for the different test cases\label{sec:numres:listofparameters}}
\end{figure}

The choice of the parameter subset size and the number of time steps are indicated in Figure \ref{sec:numres:listofparameters}. We briefly outline here the motivation for these different decisions. As far as the temporal discretization is concerned, the calculation converges after ten time steps for all test cases considered. For a more complete visualization and analysis of the results for the solution reproduction problem, we have decided to use a grid twice as fine as in the parametric case. As for the parameter grid, we have opted to start from a 2d Cartesian grid of parameters. Our case is such that the greedy algorithm converges in less than ten iterations (see the following section). Therefore, we have chosen to consider about twenty parameters in each direction of the grid.

We consider a three-dimensional quadratic tetrahedral mesh for our numerical investigations. We provide the mesh information in the Table \ref{sec:numres:table:meshinfo}.

\begin{table}[h!]
\begin{center}
\begin{tabular}{|c|c|c|c|c|}
\hline
$\mathcal{N}_e$ & $\mathcal{N}^{\rm no}$ & $\mathcal{N}^{\rm qd}$ & $\mathcal{N}$ & $\mathcal{N}_g $ \\ \hline
11 981 & 18 446 & 59 905 & 55 338 & 359 430 \\ \hline
\end{tabular}
\end{center}
\caption{Mesh information: number of three-dimensional cells (
$\mathcal{N}_e$), number of nodes ($\mathcal{N}^{\rm no}$), number of three-dimensional quadrature points ($\mathcal{N}^{\rm qd}$), size of the discretized displacement ($\mathcal{N}$) and stress vectors ($\mathcal{N}_g$ )\label{sec:numres:table:meshinfo}}
\end{table}

\subsection{Solution reproduction problem\label{sec:numres:subsec:srpb}}

We first present numerical results for a fixed configuration of parameters, namely for the centroid $\overline{\mu}\in \mathcal{P}$ to validate our ROM strategy.

\begin{figure}[h!]
\begin{center}
\centering
\begin{subfigure}[b]{0.35\textwidth}
\includegraphics[scale=1]{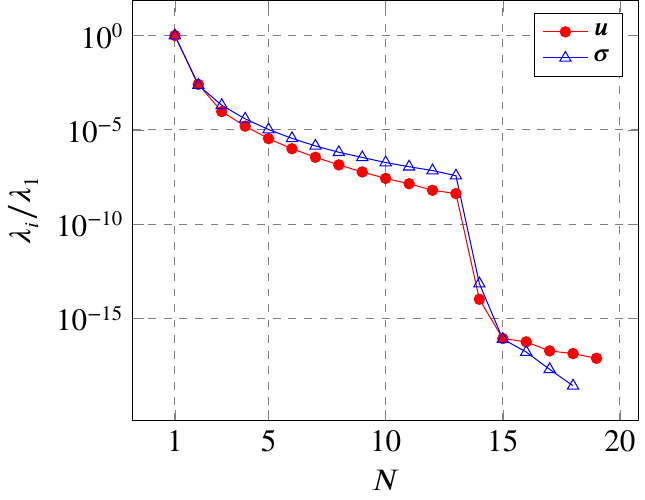} 
\caption{POD eigenvalues\label{sec:numres:subsec:srpb:fig:eigenvalues}}
\end{subfigure}
\qquad
\centering
\begin{subfigure}[b]{0.35\textwidth}
\includegraphics[scale=1]{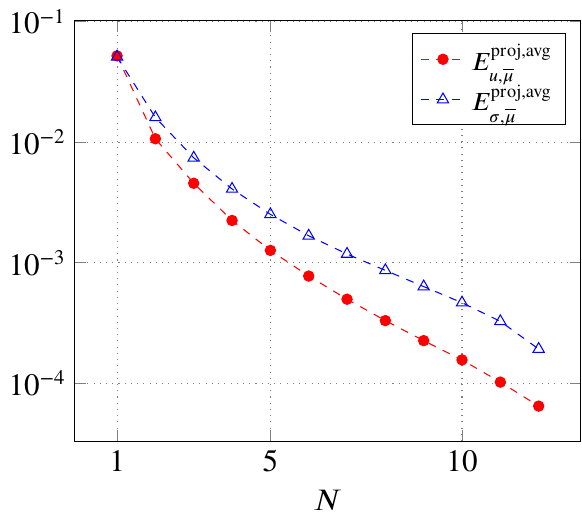} 
\caption{Projection errors\label{sec:numres:subsec:srpb:fig:projection}}
\end{subfigure}
\end{center}
\caption{Solution reproduction problem: (\subref{sec:numres:subsec:srpb:fig:eigenvalues}) behavior of the POD eigenvalues for displacement ($\mathbf{u}$) and stress ($\bm{\sigma}$) for several values of $N$($N=N_u$ for $\mathbf{u}$ and $N=N_\sigma$ for $\bm{\sigma}$; \subref{sec:numres:subsec:srpb:fig:projection}) behavior of the average projection errors) (cf. Eq.\eqref{sec:numerical:eq:projerrors})
\label{sec:numres:subsec:srpb:fig:projectionANDeigenvalues}}
\end{figure}

Figure \ref{sec:numres:subsec:srpb:fig:projectionANDeigenvalues} represents the eigenvalues obtained for the displacements and the stresses. We notice that the decays of the eigenvalues have a similar profile, although the decay of the eigenvalues is slightly faster for the displacement field than for the stress field. The plot of the projection errors as a function of the number of modes used to build the reduced space highlights this capacity to better estimate the displacement trajectory for a smaller number of modes. This suggests that in order to get projection errors in displacements and stresses at a given order of magnitude, it is mandatory to have more stress modes than displacement modes. 

\begin{figure}[h!]
\begin{center}
\includegraphics[scale=1]{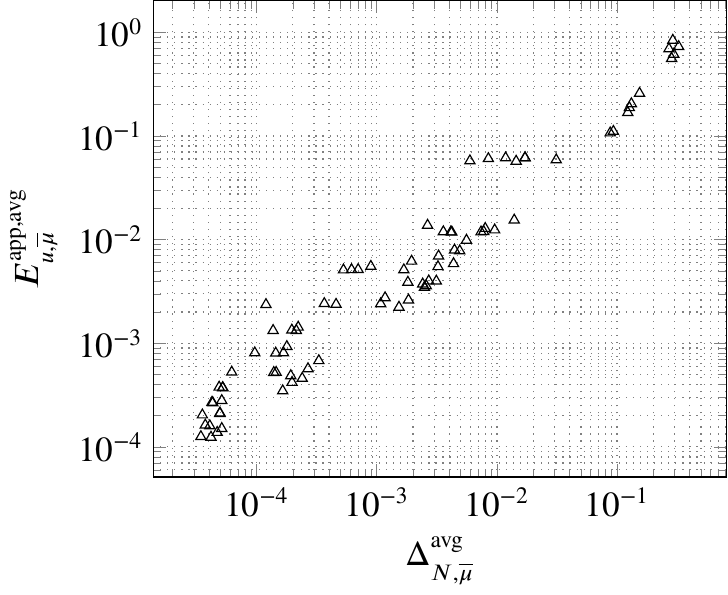} 
\end{center}
\caption{Correlation between the approximation error on the displacement ($E_{u, \overline{\mu}}^{\rm app, avg}$) and the error indicator ($\Delta_{N,\overline{\mu}}^{\rm avg}$) for the solution reproduction problem\label{sec:numres:subsec:srpb:fig:correlationapperrerrind}}
\end{figure}

Figure \eqref{sec:numres:subsec:srpb:fig:correlationapperrerrind} displays a good correlation between the the error indicator used and the approximation error on the solution fields. We point out that, in every case reported here, we have chosen to deal with all the available stress modes. Indeed, for extremly underresolved reduced spaces, the error indicator is found to be inaccurate. Since the construction of our error indicator relies on an approximation of the dual norm using the decomposition of the stress field on the space of stress modes, the correlation between the error indicator and the approximation errors may be slightly degraded for too coarse approximation spaces. This choice of treating all the stress modes does not raise overfitting problems during the Gappy-POD since we have a number of modes lower than the number of elements selected during the hyper-reduction procedure, in our quite simple case. From a practical standpoint, this choice allows us not to have to play with the ratio between the two compression tolerances ($\varepsilon_u$ and $\varepsilon_\sigma$) for the construction of the reduced problem.

\begin{figure}[h!]
\begin{center}
\includegraphics[scale=1]{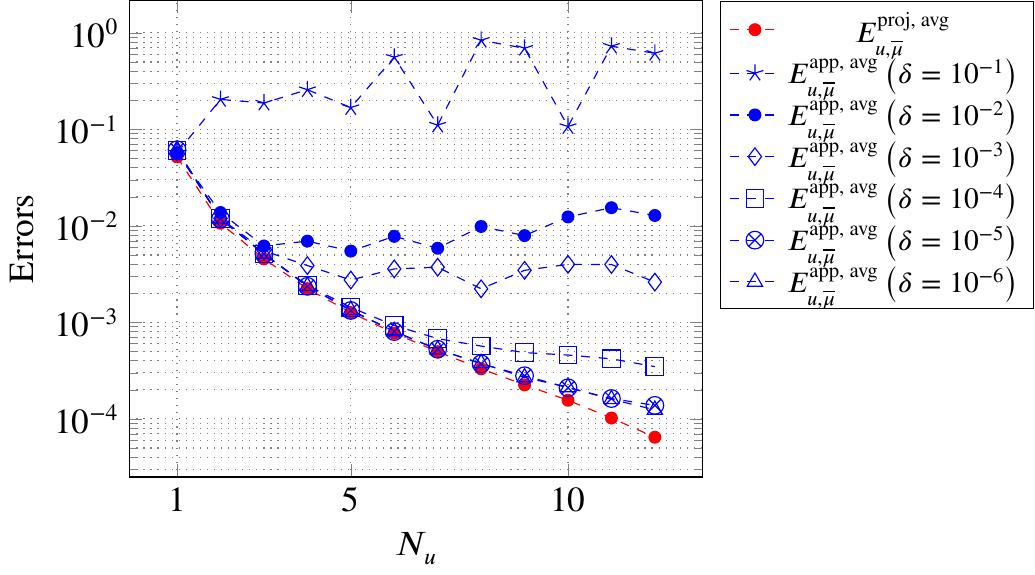} 
\end{center}
\caption{Comparison between the approximation and the projection errors with respect to the number of displacement modes $N_u$ for the solution reproduction problem. Approximation errors provided here have been computed for different values of the hyperreduction parameter $\delta$ (EQ tolerance) \label{sec:numres:subsec:srpb:fig:correlationappprojerrors}}
\end{figure}

\begin{figure}[h]
\begin{center}
\centering
\begin{subfigure}[b]{0.4\textwidth}
\includegraphics[scale=1]{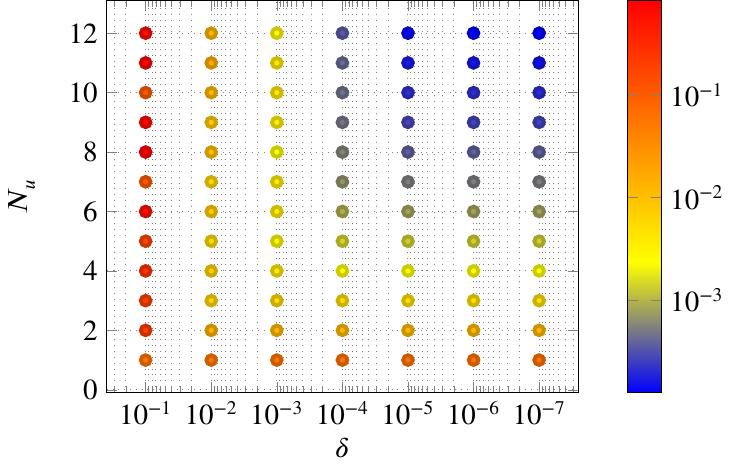} 
\caption{$E_{u, \overline{\mu}}^{\rm app, avg}$\label{sec:numres:subfigure:approxerror}}
\end{subfigure}
\qquad
\centering
\begin{subfigure}[b]{0.4\textwidth}
\includegraphics[scale=1]{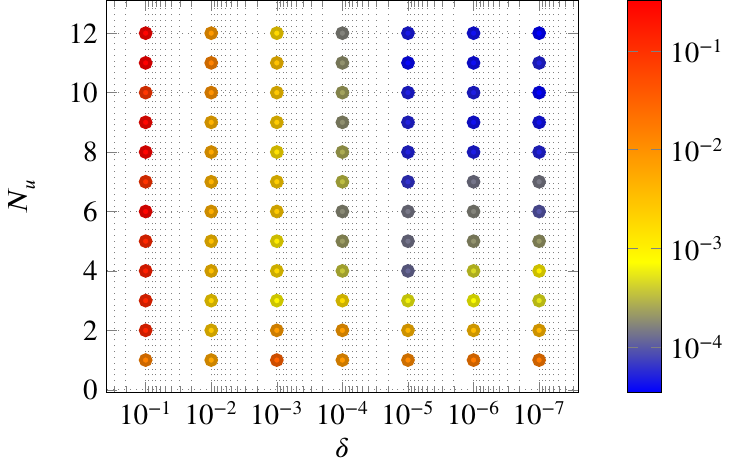} 
\caption{$\Delta_{N,\mu}^{\rm avg}$\label{sec:numres:subfigure:errind}}
\end{subfigure}
\end{center}
\begin{center}
\centering
\begin{subfigure}[b]{0.4\textwidth}
\includegraphics[scale=1]{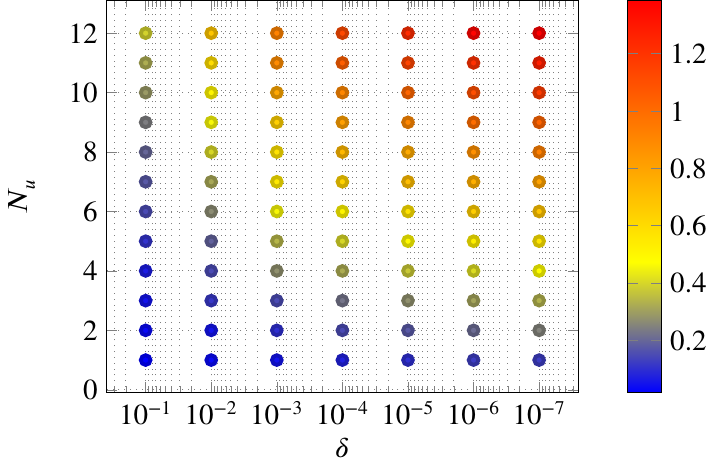} 
\caption{$\norm{\bm{\rho}^{\rm eq }}_{\ell_0} / \ \norm{\bm{\rho}^{\rm hf}}_{\ell_0} \ [\%]$\label{sec:numres:subfigure:pse}}
\end{subfigure}
\qquad
\centering
\begin{subfigure}[b]{0.4\textwidth}
\includegraphics[scale=1]{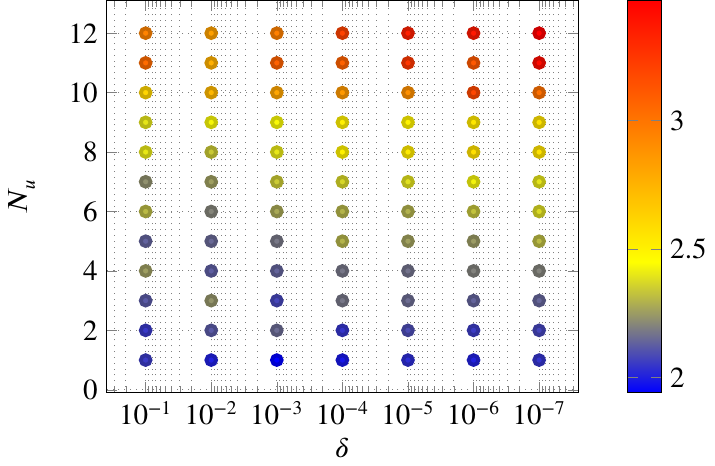}
\caption{ROM CPU cost / HF CPU cost $[\%]$\label{sec:numres:subfigure:pcpu}}
\end{subfigure}
\end{center}
\caption{Solution reproduction problem: colormaps of: (\subref{sec:numres:subfigure:approxerror}) approximation errors, (\subref{sec:numres:subfigure:errind}) error indicators , (\subref{sec:numres:subfigure:pse}) percentage of selected elements, (\subref{sec:numres:subfigure:pcpu}) percentage of CPU time for different size of reduced order basis and hyper-reduction parameters \label{sec:numres:fig:allcolormaps}}
\end{figure}

We have built reduced models for various numbers of modes (compression of the solution space) and various hyper-reduction parameters (size of the reduced mesh). The aim of investigating this grid of hyperparameters of the reduced model is multifaceted. First, it enables to investigate a wide range of approximation errors. Indeed, the quality of the approximated solution depends on the approximation quality of the integrals involved in the problem ($\delta$) but also on the approximation quality of the trajectory ($N_u$). This variation allows us to highlight the correlation between the approximation error on the displacement field and the error indicator that we have presented (Figure \ref{sec:numres:subsec:srpb:fig:correlationapperrerrind}, and Colormaps \ref{sec:numres:subfigure:approxerror} and \ref{sec:numres:subfigure:errind}). Moreover, for a fixed number of modes, the projection error constitutes a theoretical lower bound that we wish to be able to reach by solving the reduced problem. However, the hyper-reduction process introduces a new approximation. In Figure \ref{sec:numres:subsec:srpb:fig:correlationappprojerrors}, we illustrate that the approximation error tends towards the projection error for small $\delta$ values, while a less restrictive parameter degrades the solution ($\delta = 10^{-1}$ for example). The slight differences between approximation and projection errors between the last two $\delta$ values comes from the fact that we hit the tolerance of the iterative Newton algorithm used in the HF solver (which is chosen as the same as in the reduced solver).

\begin{figure}[h!]
\begin{center}
\includegraphics[scale=1]{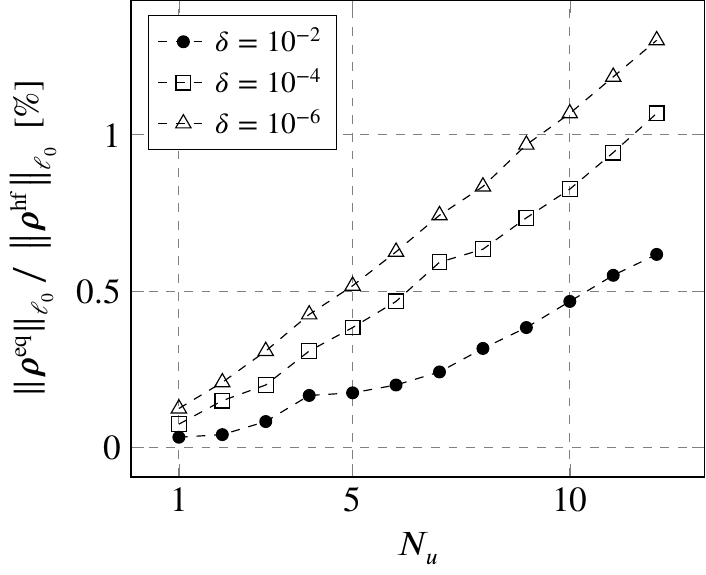} 
\end{center}
\caption{Percentage of three-dimensional selected elements depending on the size of the primal basis ($N_u$) and the EQ tolerance $\delta$ for the solution reproduction problem}
\end{figure}

Much more, we observe that the empirical quadrature procedure is able to significantly reduce the size of the mesh used for online calculations. We keep at most a few percent of the number of elements in the HF mesh. We thus drastically reduce the computational cost compared to a HF problem. Indeed, the cost of a reduced problem represents only a few percent of the cost of the HF calculation. The computational cost reduction, correlated to the number of selected elements (Colormap \ref{sec:numres:subfigure:pse}), depends both on the number of selected modes and on the hyper-reduction parameter that we choose.

\begin{remark}
It should be noticed in the following case that the ratio between the computational cost of the reduced problem and the percentage of selected elements are not strictly correlated, even if the two quantities follow the same tendency. Indeed, from an algorithmic point of view, the reduction of the mesh is not the only operation involved between a HF computation and a reduced computation, since the projection on the modes and the hyper-reduction entail a modification of the size of the system, but also of the conditioning of the latter (this can lead to more Newton iterations for a reduced computation for example). Furthermore, the implementation has been done in an industrial code where fixed costs related to verification and memory allocation processes are necessary whatever the computation. Nevertheless, in Figure \ref{sec:numres:fig:allcolormaps}, we provide a numerical validation that the percentage of selected elements gives us a good hint on the gain in terms of computational cost.
\end{remark}

We have observed that the approximation error on the stresses follows the same pattern as the approximation error on the displacements on the hyperparameter grid ($N_u \times \delta$). This comment brings us to report only approximation errors on displacements in this contribution.

\subsection{Parametric problem}

\subsubsection{Parametric problem $\mu=\nu$\label{sec:numres:subsec:param:subsubsec:nu}}

In this section, we consider the variation of a single parameter ($\nu$), for a training set of $\left|\Theta_{{\rm train}, \nu}\right| = 20$ values of this parameter. The numerical results presented here and in the last sub-section were performed for a smaller number of time steps for the sake of efficiency (see parameters in Table \ref{sec:numres:listofparameters}). 

\begin{figure}[h!]
\begin{center}
\includegraphics[scale=1]{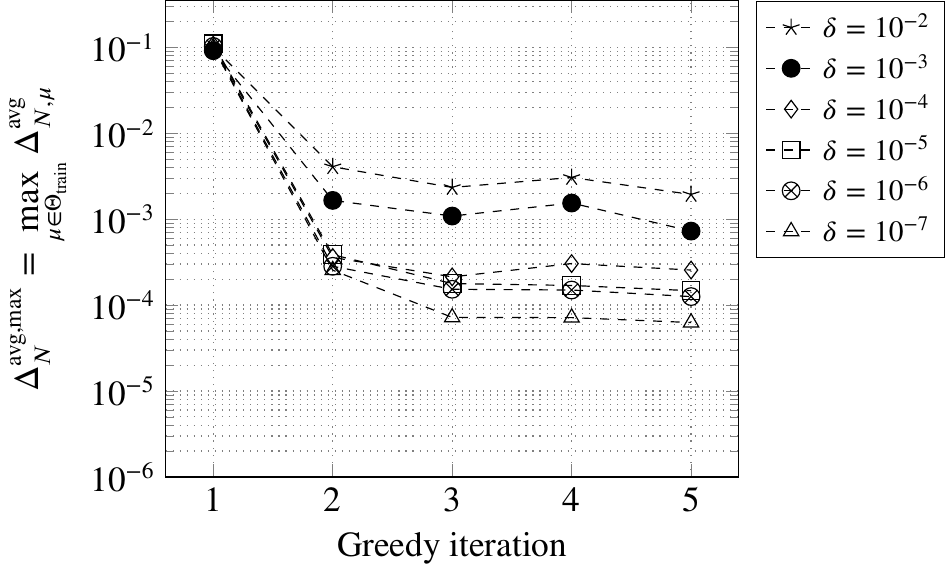} 
\end{center}
\caption{Maximum of the time-averaged error indicator over the training set depending on the Greedy iteration for the parametric problem $\mu = \nu$\label{sec:numres:subsec:param:subsubsec:nu:figure:maxerrind}} 
\end{figure}

We were able to test the greedy approach on this single parameter nonlinear case where we fixed a given number of iterations ($N_{\rm Gr} = 5$). For the studied example, the algorithm has reached its convergence for the following number of iteration. All the examples reported here have been carried out for the tolerance $\varepsilon_u = 10^{-5}$, which ensures a good approximation error on the explored parameters. The evolution of the maximum of the error indicator (Figure \ref{sec:numres:subsec:param:subsubsec:nu:figure:maxerrind}) for several hyper-reduction parameters shows a convergence after a few iterations. The plateau reached by the error indicator differs with the accuracy of the approximation of the integrals of the problem.

\begin{figure}[h!]
\begin{center}
\centering
\begin{subfigure}[b]{0.4\textwidth}
\includegraphics[scale=1]{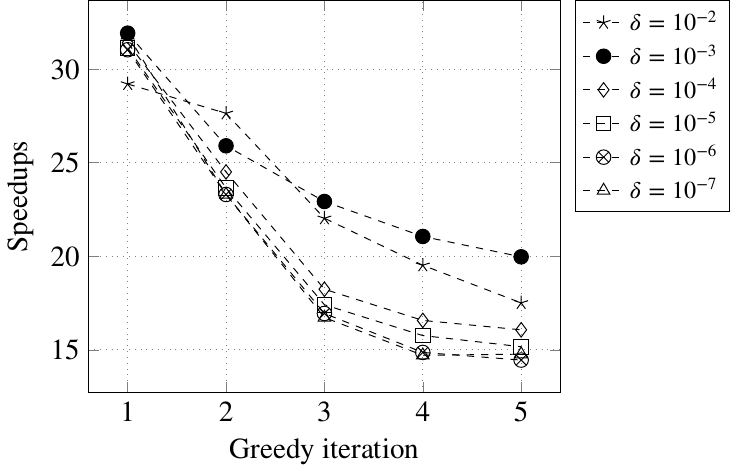} 
\end{subfigure}
\qquad
\centering
\begin{subfigure}[b]{0.4\textwidth}
\includegraphics[scale=1]{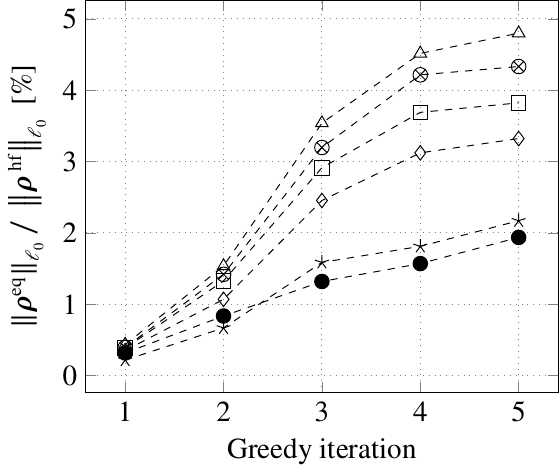} 
\end{subfigure}
\end{center}
\begin{center}
\centering
\begin{subfigure}[b]{\textwidth}
\begin{center}
\includegraphics[scale=1]{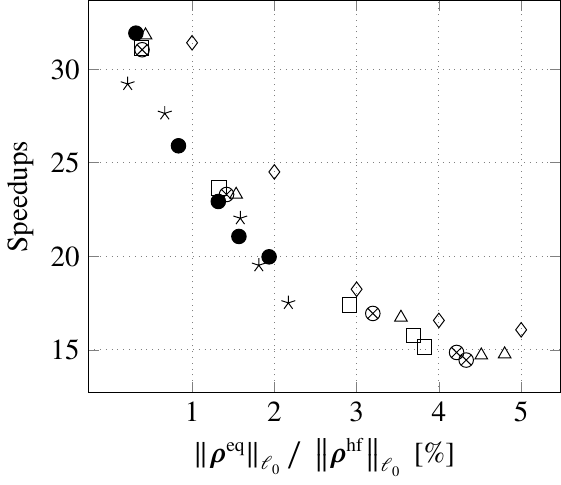} 
\end{center}
\end{subfigure}
\end{center}
\caption{Parametric problem $\mu = \nu$: informations (speedups and percentage of three-dimensional selected elements $\norm{\bm{\rho}^{\rm eq}}_{\ell_0} / \ \norm{\bm{\rho}^{\rm hf}}_{\ell_0} \ [\%]$) through the Greedy iterations ($\varepsilon=10^{-5}$) for different values of $\delta$\label{sec:numres:subsec:param:subsubsec:nu:figure:infogreedy}}
\end{figure}

We report the computational costs associated with the reduced solver by giving the speedup=$\frac{\rm HF cost}{\rm ROM cost}$, where the HF cost is the computational time of solving the HF problem whereas rom cost is the online cost of evaluating the problem. In Figure \ref{sec:numres:subsec:param:subsubsec:nu:figure:infogreedy}, we notice that the gain in computation time decreases with each iteration as the percentage of selected elements increases with the number of HF problems to be estimated and the number of modes to be included in the reduced basis. Nevertheless, for this single-parameter problem, the speedups obtained are always higher than 10 or even 15, which implies a drastic decrease of the computation time for the model evaluation. Moreover, the parametric manifold is in our case very well approximated after a small number of iterations. The plateau observed in Figure \ref{sec:numres:subsec:param:subsubsec:nu:figure:maxerrind} is reached after a few iterations and shows that for the given tolerance of hyperreduction and the desired precision in the compression of the base, there is no more gain in exploring a new parameter.

\begin{figure}[h!]
%\centering
\begin{subfigure}[c]{0.25\textwidth}
\begin{center}
\includegraphics[scale=0.3]{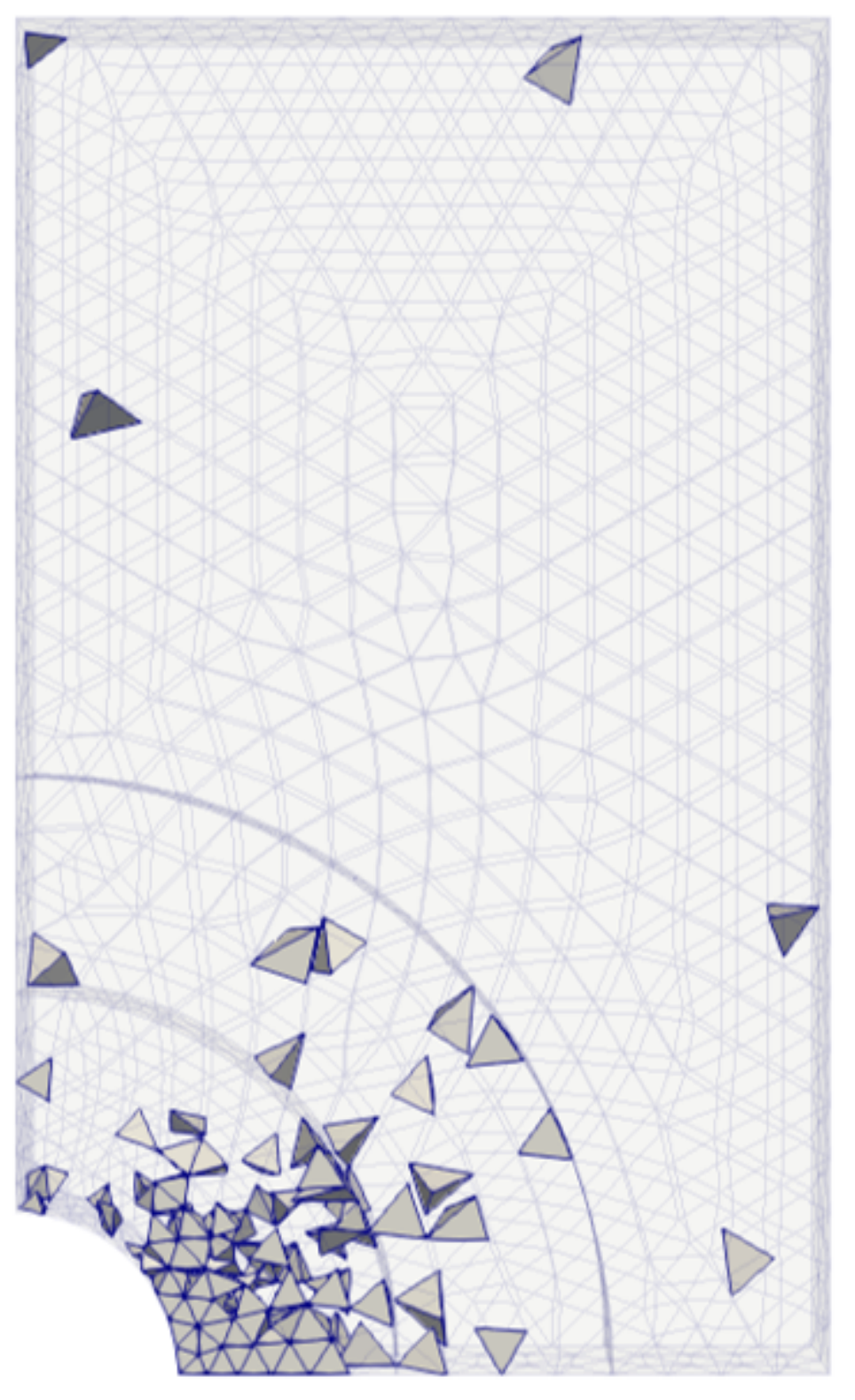} 
\end{center}
\caption{Reduced mesh for $\delta=10^{-2}$}
\end{subfigure}
\quad 
\begin{subfigure}[c]{0.25\textwidth}
\begin{center}
\includegraphics[scale=0.3]{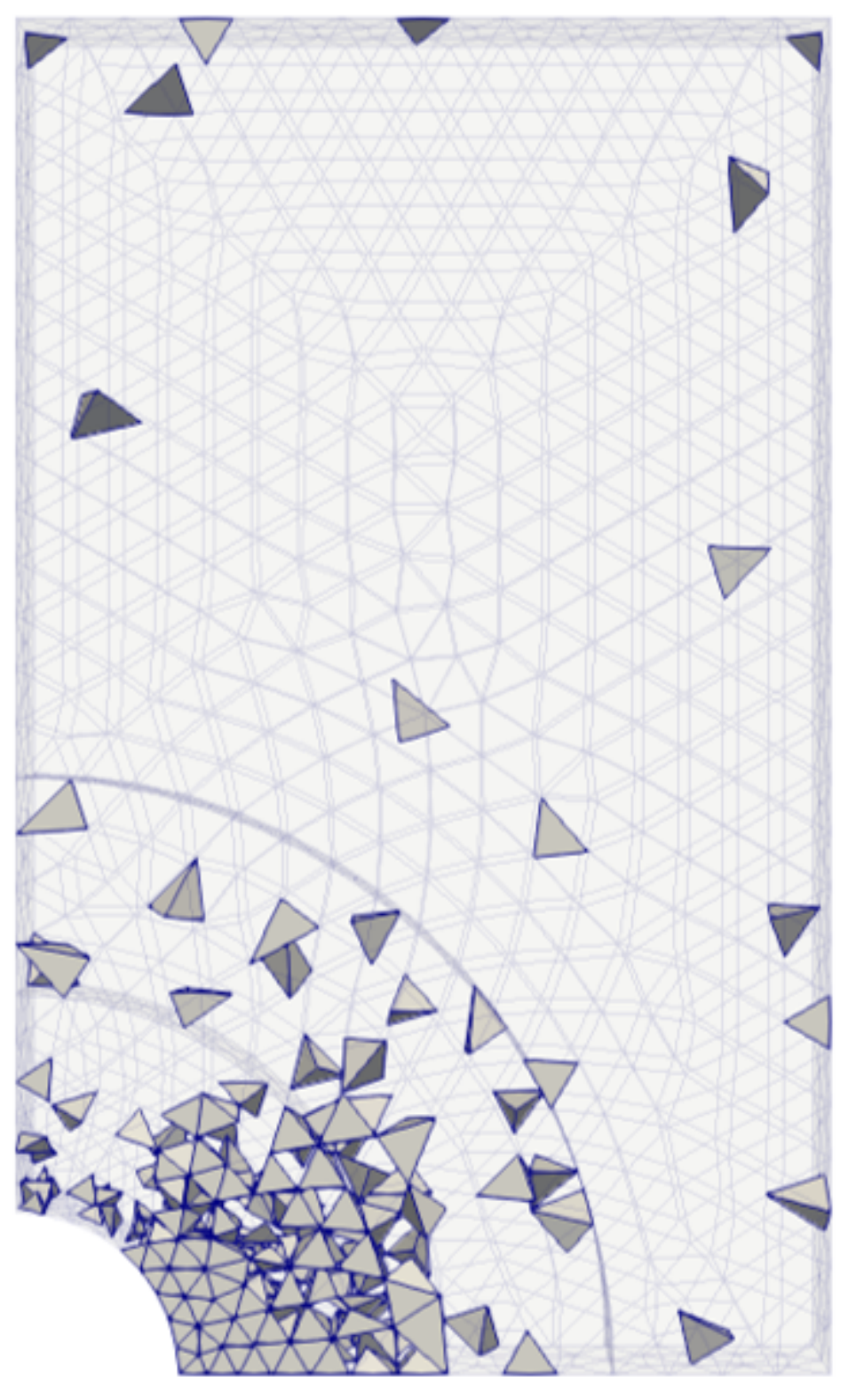} 
\end{center}
\caption{Reduced mesh for $\delta=10^{-7}$}
\end{subfigure}
\qquad 
\begin{subfigure}[c]{0.3\textwidth}
\begin{center}
\begin{tabular}{|c|c|c|}
\hline
$\delta$ & $10^{-2}$ & $10^{-7}$ \\ \hline
Speedup & 17.51 &  14.77  \\ \hline
$\norm{\bm{\rho}^{\rm eq}}_{\ell_0} / \ \norm{\bm{\rho}^{\rm hf}}_{\ell_0} \ [\%]$ & 2.17 & 4.80 \\ \hline
$N_u$ & 23 & 22 \\ \hline
\end{tabular}
\caption{Summary of outputs of the hyper-reduction procedure}
\end{center}
\end{subfigure}
\caption{Parametric problem $\mu = \nu$: hyper-reduced ROM and selected elements at the end of the POD-Greedy procedure ($\varepsilon_u=10^{-5}$) for two different values of $\delta$\label{sec:numres:subsec:param:subsubsec:nu:figure:redmesh}}
\end{figure}

We notice on the reduced meshes obtained at the end of the numerical procedure (Figure \ref{sec:numres:subsec:param:subsubsec:nu:figure:redmesh}) that the selected elements are mainly located around the hole, which matches the region where the material enters a nonlinear regime (plastic regime). As we would expect, few elements are selected in the areas where the behavior is purely elastic (linear).

\subsubsection{Multi-parametric problem $\mu = \left(\nu, a_{\rm pui}\right)$\label{sec:numres:subsec:param:subsubsec:nuapui}}

\begin{figure}[h]
\begin{center}
\centering
\begin{subfigure}[b]{0.4\textwidth}
\includegraphics[scale=1]{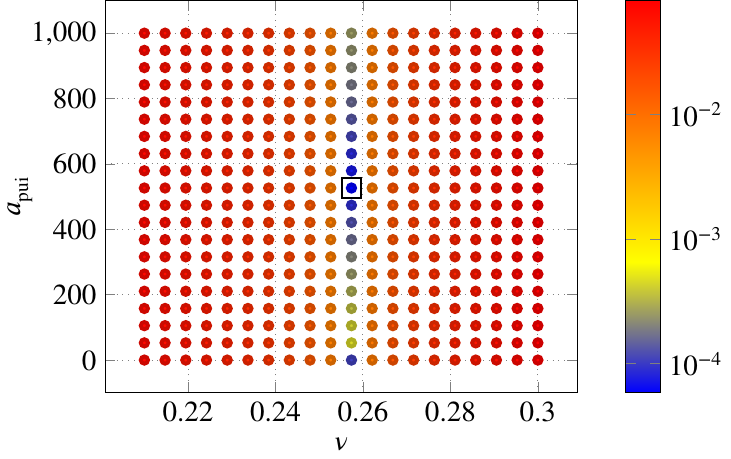}
\caption{Greedy iteration n$^\circ$1}
\end{subfigure}
\qquad
\centering
\begin{subfigure}[b]{0.4\textwidth}
\includegraphics[scale=1]{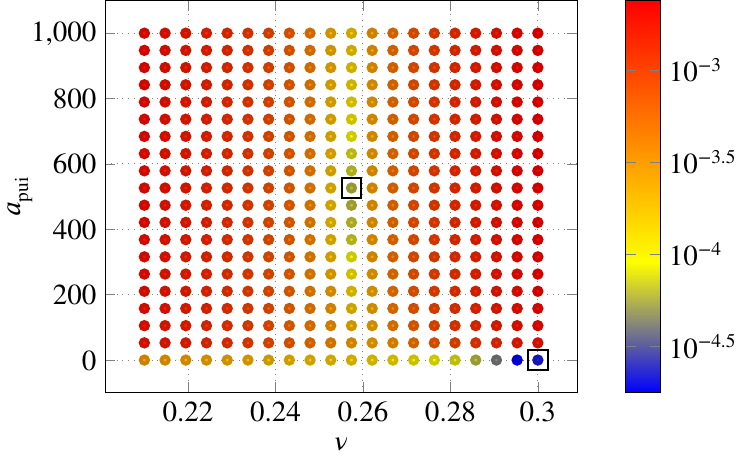} 
\caption{Greedy iteration n$^\circ$2}
\end{subfigure}
\end{center}
\begin{center}
\centering
\begin{subfigure}[b]{0.4\textwidth}
\includegraphics[scale=1]{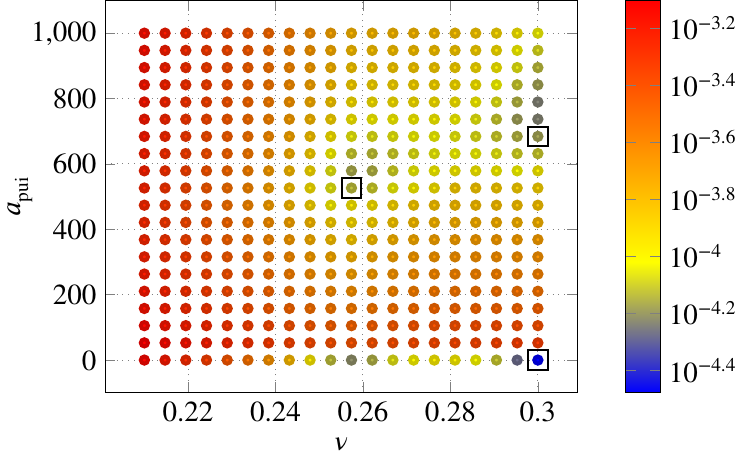} 
\caption{Greedy iteration n$^\circ$3}
\end{subfigure}
\qquad
\centering
\begin{subfigure}[b]{0.4\textwidth}
\includegraphics[scale=1]{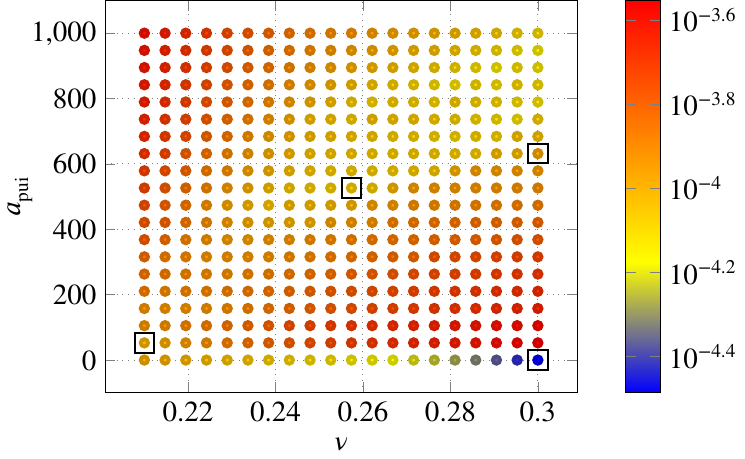} 
\caption{Greedy iteration n$^\circ$4}
\end{subfigure}
\end{center}
\caption{Parametric problem $\mu = \left(\nu, a_{\rm pui}\right)$: colormaps of the time-averaged error indicators and selected parameters (points squared in black) for every Greedy iterations\label{sec:numres:subsec:param:subsubsec:nuapui:errindevol}}
\end{figure}

\begin{figure}[h]
\begin{center}
\centering
\begin{subfigure}[c]{0.28\textwidth}
\includegraphics[scale=1]{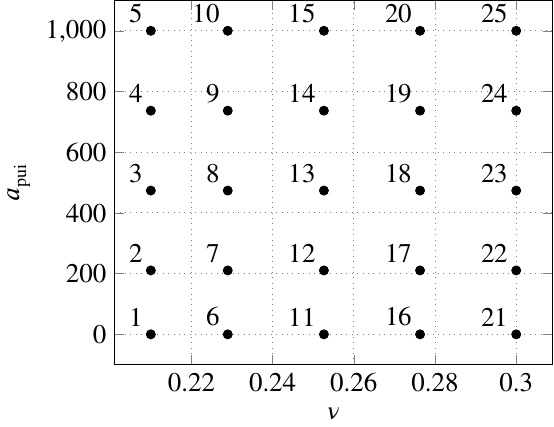} 
\caption{Test grid for error computation and indexing of the test set\label{fig:nuapui:errindapprox:repartition}}
\end{subfigure}
\qquad
\centering
\begin{subfigure}[c]{0.28\textwidth}
\includegraphics[scale=1]{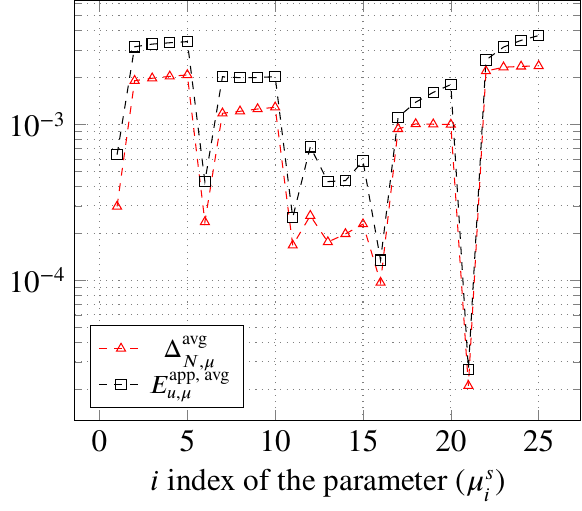} 
\caption{Greedy iteration n$^\circ$2}
\end{subfigure}
\qquad
\centering
\begin{subfigure}[c]{0.28\textwidth}
\includegraphics[scale=1]{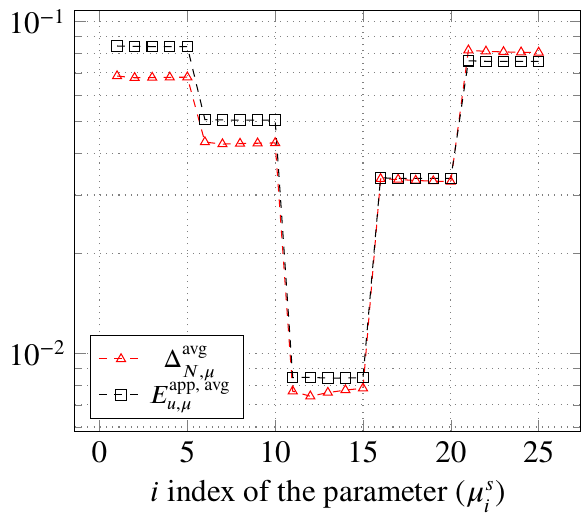} 
\caption{Greedy iteration n$^\circ$1}
\end{subfigure}
\end{center}
\caption{Parametric problem $\mu = \left(\nu, a_{\rm pui}\right)$: correlation between the error indicator and the approximation errors depending on the numerotation of the test parameters (repartition shown in \ref{fig:nuapui:errindapprox:repartition}) \label{sec:numres:subsec:param:subsubsec:nuapui:errindapprox}}
\end{figure}

Finally, we provide a numerical example for two parameters. We address a training set of size $\left|\Theta_{\rm train}\right| = 400$. We report here the results for tolerances $\varepsilon_u = 10^{-5}$ and $\delta=10^{-7}$. The choice of the hyper-reduction parameter is chosen here of the same order of magnitude as the Newton tolerance for HF computation. In Figure \ref{sec:numres:subsec:param:subsubsec:nuapui:errindevol}, we present the evolution of the error indicator we compute over the greedy iterations.  By comparing, for instance, the colormap at the second iteration and at the third iteration, we notice that the sampling of a parameter leads to a decrease of the indicator value in the neighborhood of the given parameter. Moreover, we were interested in the correlation of the error indicator with the error indicator especially with out-of-sample parameters. To this end, we defined a sub-grid of 25 points, (5$\times$5 Cartesian grid of the parameters), on which the HF calculations were performed in order to dispose of the projection error. In Figure \ref{sec:numres:subsec:param:subsubsec:nuapui:errindapprox}, we show the profiles obtained for the error indicators and for the approximation errors on the parameters chosen for the test. It appears that the error indicator seems to follow the behaviour of the approximation error.

We point out that we only provide results on the first iterations because we have limited ourselves to a small number of iterations as the problem is sufficiently well approximated in a short time. It would therefore not be relevant to compare relative errors where the variation between the parameters becomes insignificant.

\section{Conclusion}\label{sec:conclusion}

In this work, we developed and numerically validated a model order reduction procedure for a generic class of nonlinear mechanical problems with internal variables. We successfully implemented the method directly with an HF industrial code and validated it on an elastoplastic material. We proposed a time-averaged error indicator to drive the offline Greedy sampling, which is cost-efficient and has been shown numerically correlated to the approximation errors, and  we developped an element-wise empirical quadrature procedure to reduce online costs. The whole procedure delivered impressive computational cost improvements cost improvements in the order of $O(20-25)$ with relative prediction errors in the order of $10^{-3}$.

We aim to extend our methodology in several directions. First, we try to extend this approach to more complex problems with more marked differences depending on the physical parameters used. Indeed, the algorithm presented here allows to approach the parametric variety after only a few iterations. We therefore seek to test the approach on more complex problems to assess the relevance of our approach for other nonlinear quasi-static mechanical problems and we wish highlight the interest of this adaptive approach compared to naive approach (cartesian grid) provided that we have a more complex physical problem to address. We also wish to explore more sophisticated sampling strategies to reduce offline costs. Moreover, we are looking to extend this approach to real-world challenging industrial problems, such as with three-dimensional and one-dimensional mechanical couplings, namely in the case of engineering studies involving prestressed concrete.

\section*{Acknowledgement}

This work was partly funded by ANRT (French National Association for Research
and Technology) and EDF. We wish to express our thanks to the \textsf{code$\_$aster} development team and the contributors to this code. Our efforts are related to the use and development of the python library Mordicus\cite{mordicus} (funded by a 'French Fonds Unique Interministériel' FUI project) which is intended to provide a tool for the development of model reduction methods for industrial applications\cite{grosjean2022variations}\cite{daniel2021machine}.

\clearpage
\newpage

%\nocite{*}% Show all bib entries - both cited and uncited; comment this line to view only cited bib entries;

\appendix

\section{Newton solver}\label{app:newton}

In this appendix, we provide details on the numerical procedure used for solving nonlinear systems with dualisation of boundary conditions in the \textsf{code$\_$aster} framework. For this purpose, we first discuss the procedure used when the kinematic conditions are handled by Dirichlet elimination, before introducing the dualization of the boundary conditions and the stopping criteria considered.

\subsection{No dualisation of the boundary conditions}

We focus on looking for the $k$-th timestep solution . The resolution is performed by a Newton-Raphson type algorithm, which is an incremental algorithm. The iterative process is driven by the search for a solution at each iteration according to the knowledge at the previous iteration:

$$\mathbf{u}^{(k)}_{\theta + 1} = \mathbf{u}^{(k)}_{\theta}+\Delta \mathbf{u}^{(k)}_{\theta}$$

The iterate is computed from the solution of the linear system, expressed with the Jacobian matrix (also called tangent matrix in mechanics) evaluated in $\mathbf{u}^{(k)}_{\theta}$:

$$\mathbf{R}^{\rm hf}\left(\mathbf{u}^{(k)}_{\theta + 1}\right) \approx \mathbf{R}^{\rm hf}\left(\mathbf{u}^{(k)}_{\theta}\right) + \mathbf{K}^{(k)}_{\theta} \cdot \Delta \mathbf{u}^{(k)}_{\theta}  = 0, \quad \text{with} \quad \mathbf{K}^{(k)}_{\theta} = \frac{\mathbf{R}^{\rm hf}}{\partial \mathbf{u}} \left(\mathbf{u}^{(k)}_{\theta}\right)$$

\subsection{Dualization of the boundary conditions}

For the dualisation of constraints, we must investigate a new increment in displacement and in terms of Lagrange multipliers:

\begin{equation}
\left\{
\begin{array}{rcl}
\mathbf{u}^{(k)}_{\theta + 1} &=& \mathbf{u}^{(k)}_{\theta}+\Delta \mathbf{u}^{(k)}_{\theta}\\
\bm{\lambda}^{(k)}_{\theta + 1}&=& \bm{\lambda}^{(k)}_{\theta } + \Delta \bm{\lambda}^{(k)}_{\theta}
\end{array}
\right.
\end{equation}

The task is hence to solve the following nonlinear system:

\begin{equation}
\left\{
\begin{array}{rcl}
\mathbf{R}^{\rm hf}\left(\mathbf{u}^{(k)}_{\theta + 1}\right)  + \mathbf{B}^T \bm{\lambda}^{(k)}_{\theta + 1}&=& 0\\
\mathbf{B}\mathbf{u}^{(k)}_{\theta + 1}&=& \mathbf{u}_d^{(k)}
\end{array}
\right.
\end{equation}

Using a linearization analogous to the equation, and exploit the linearity of the operator associated with the kinematic conditions, the assembled discretized system (for one $\theta$ iteration) is decomposed as:

\begin{equation}
\left\{
\begin{array}{rcl}
\mathbf{R}^{\rm hf}\left(\mathbf{u}^{(k)}_{\theta}\right) + \mathbf{K}^{(k)}_{\theta} \cdot \Delta \mathbf{u}^{(k)}_{\theta} + \mathbf{B}^T \bm{\lambda}^{(k)}_{\theta} + \mathbf{B}^T \Delta \bm{\lambda}^{(k)}_{\theta}&=& 0\\
\mathbf{B}\mathbf{u}^{(k)}_{\theta} + \mathbf{B}\Delta\mathbf{u}^{(k)}_{\theta}&=& \mathbf{u}_d^{(k)}
\end{array}
\right.
\end{equation}

\noindent which leads to the following saddle-point problem:

\begin{equation}
\begin{bmatrix}
\mathbf{K}^{(k)}_{\theta} & \mathbf{B}^T \\
\mathbf{B} & 0 
\end{bmatrix}
\begin{bmatrix}
 \Delta \mathbf{u}^{(k)}_{\theta}\\
\Delta \bm{\lambda}^{(k)}_{\theta}
\end{bmatrix} =
\begin{bmatrix}
-\mathbf{R}^{\rm hf}\left(\mathbf{u}^{(k)}_{\theta}\right) -\mathbf{B}^T \bm{\lambda}^{(k)}_{\theta}  \\
 \mathbf{u}_d^{(k)} - \mathbf{B}\mathbf{u}^{(k)}_{\theta}
\end{bmatrix}
\end{equation}

\subsection{Stopping criterium}
\label{app:newton:subsec:stopcrit}

Under the philosophy of the formulations in \textsf{code$\_$aster}, the internal contributions (work of internal forces) and external contributions (forces applied to the system) are evaluated separately in the assembled residue:

\begin{equation}
\mathbf{R}^{\rm hf}\left(\mathbf{u}^{(k)}\right) = \mathbf{F}^{{\rm int}, (k)}\left(\mathbf{u}^{(k)}\right) - \mathbf{F}^{{\rm ext}, (k)}
\end{equation}

Different criteria are available in \textsf{code$\_$aster}. The reader may refer to the code documentation for more details. Our choice is a relative criterion defined as follows:

\begin{equation}
\frac{\norm{\mathbf{R}^{\rm hf}\left(\mathbf{u}^{(k)}_{\theta}\right)  + \mathbf{B}^T \bm{\lambda}^{(k)}_{\theta}}_\infty}{\norm{ \mathbf{B}^T \bm{\lambda}^{(k)}_{\theta}- \mathbf{F}^{{\rm ext}, (k)}}_\infty} \leq \varepsilon_{\rm newt}
\end{equation}

The vector $\mathbf{B}^T \bm{\lambda}^{(k)}_{\theta}$ can be interpreted physically as the opposite of the support reactions at the nodes where the conditions are dualised. The convergence criterion can be seen as a process of normalizing the residual calculated at a given iteration with respect to the forces exerted on the system at that iteration (external forces and support reactions).

\section{Dictionary construction}\label{app:eqconstruction}

\subsection{Solution reproduction problem example}

\subsubsection{Formulation}

We resume the example introduced in the section \ref{sec:metho:subsec:srpb}, i.e. the case of a solution reproduction problem. We describe more precisely the hyper-reduction process used in our methodology. We keep the same notations as previously introduced. In such a scenario, we have $K$ HF snapshot (displacements and stresses) and $N_u$ primal modes at our disposal. We hence have $n_{\rm int} = K\times N_u$ \emph{manifold accuracy constraints} to fulfill:

\begin{equation}
\left(\mathbf{G}\right)_{{\rm lines}(n, k) , \ q} = \mathcal{R}^{\sigma, \rm hf}_{q}\left(\mathbf{E}_q^{\rm qd}\bm{\sigma}^{(k)}, \ \mathbf{E}_q^{\rm no}\bm{\zeta}_{u,n}\right) \quad \text{and} \quad \left(\mathbf{y}\right)_{{\rm lines}(n, k)} = \mathcal{R}^{\sigma, \rm hf}\left( \mathbf{E}_q^{\rm qd}\bm{\sigma}^{(k)}, \ \mathbf{E}_q^{\rm no}\bm{\zeta}_{u,n}\right)
\end{equation}

\noindent where $\mathbf{G}\in \mathbb{R}^{n_{\rm int} \times N_e}$ and $\mathbf{y}\in \mathbb{R}^{n_{\rm int}} $ and $\rm{lines} : \ (k,n)\in \mathbb{R}^{K \times N_u} \rightarrow \mathbb{R}^{n_{\rm int}}$ a bijection used to have a unique numerotation of rows (set by the way we build the dictionnary). The last row of the dictionnary is set in order to fulfill the \emph{constant-function constraint}:

\begin{equation}
\left(\mathbf{G}\right)_{n_{\rm int} + 1 , \ q} = \left|K_q\right|, \quad \text{and} \quad \left(\mathbf{y}\right)_{n_{\rm int} + 1 } = \left|\Omega\right|
\end{equation}

\subsubsection{Separation of integrals}

As we restrict ourselves to a single-mesh study, we have only volumic forces applied to the system. From a practical viewpoint, adding directly the residuals can load to numerical instabilities. Indeed, if the probelm is well represented by a single mode, $\zeta_{u, n*}$, we can have: 

$$\mathcal{R}^{\sigma, \rm hf}\left(\mathbf{E}_q^{\rm qd}\bm{\sigma}^{(k)}, \ \mathbf{E}_q^{\rm no}\bm{\zeta}_{u,n}\right) \approx 0 $$

To tackle this issue, we chose to split the residual in two contributions: one for the internal forces and the other for the external forces. Such an implementation is consistent with \textsf{code$\_$aster} discrete formulation. The residuals can be expressed thanks to the variationnal form as:

\begin{equation}
\mathcal{R}^{\sigma, \rm hf}_{q}\left(\mathbf{E}_q^{\rm qd}\bm{\sigma}^{(k)}, \ \mathbf{E}_q^{\rm no}\bm{\zeta}_{u,n}\right) = \int_{\Omega_q} \sigma^{(k)}:\nabla_s\zeta_{u,n}  \ dx - \int_{\Omega_q} f_v \cdot \zeta_{u,n} dx
\end{equation}

We can then define the contributions:

\begin{equation}
\left\{
\begin{array}{rcl}
\mathcal{R}^{\sigma, \rm hf, \rm int}_{q}\left(\mathbf{E}_q^{\rm qd}\bm{\sigma}^{(k)}, \ \mathbf{E}_q^{\rm no}\bm{\zeta}_{u,n}\right)&=&\int_{\Omega_q} \sigma^{(k)}:\nabla_s\zeta_{u,n}  \ dx \\
\mathcal{R}^{\sigma, \rm hf, \rm ext}_{q}\left(\mathbf{E}_q^{\rm qd}\bm{\sigma}^{(k)}, \ \mathbf{E}_q^{\rm no}\bm{\zeta}_{u,n}\right)&=&\int_{\Omega_q} f_v \cdot \zeta_{u,n} dx
\end{array}
\right.
\end{equation}

With this formulation, we have $n_{\rm lin} = (K + 1) \times N_u$ and $\mathbf{G}$ and $\mathbf{y}$ are modified accordingly:

\begin{equation}
\left(\mathbf{G}\right)_{{\rm lines}(n, k, \rm *) , \ q} = \mathcal{R}^{\sigma, \rm hf}_{q}\left(\mathbf{E}_q^{\rm qd}\bm{\sigma}^{(k)}, \ \mathbf{E}_q^{\rm no}\bm{\zeta}_{u,n}\right) \quad \text{and} \quad \left(\mathbf{y}\right)_{{\rm lines}(n, k, \rm *)} = \mathcal{R}^{\sigma, \rm hf}\left(\mathbf{E}_q^{\rm qd}\bm{\sigma}^{(k)}, \ \mathbf{E}_q^{\rm no}\bm{\zeta}_{u,n}\right)
\end{equation}

\subsubsection{Normalization}

A challenge related to the orders of magnitude arises in the optimization problem construction. Indeed, we have lines related to volume constraints while others are related to internal or external forces. Since the algorithms convergence criteria used are designed on the residuals (in the sense of optimisation, i.e. $\norm{\mathbf{G}\bm{\rho} - \mathbf{y}}_*$), it is likely that some constraints are 'overlooked' because of the differences in order of magnitude. To ensure a good behaviour of our strategy, we normalize the whole dictionary to have an addimensionalized problem:

\begin{equation}
\left(\mathbf{G}\right)_{{\rm lines}(n, k, \rm *) , \ q} = \frac{\mathcal{R}^{\sigma, \rm hf}_{q}\left(\mathbf{E}_q^{\rm qd}\bm{\sigma}^{(k)}, \ \mathbf{E}_q^{\rm no}\bm{\zeta}_{u,n}\right)}{\mathcal{R}^{\sigma, \rm hf}\left(\mathbf{E}_q^{\rm qd}\bm{\sigma}^{(k)}, \ \mathbf{E}_q^{\rm no}\bm{\zeta}_{u,n}\right)} \quad \text{and} \quad \left(\mathbf{G}\right)_{n_{\rm int} + 1 , \ q} = \frac{\left|K_q\right|}{\left|\Omega\right|}
\end{equation}

Thus, the second member consists only of a unitary vector:

\begin{equation}
\left(\mathbf{y}\right)_{{\rm lines}(n, k, \rm *)} = 1, \quad \text{and} \quad \left(\mathbf{y}\right)_{n_{\rm int} + 1 } = 1
\end{equation}

This approach is well suited to industrial codes that are not necessarily designed to have dimensionless formulations.

\section{Error indicator}
\label{app:errorind}
\subsection{Time-dependent external forces}

We consider a formulation where the external loading can vary during time. In such a situation, we have a different linear form for each timestep. We can then define:

$$\left(\psi_{N_\sigma + 1}^{\sigma, (k)}, \ v\right) = \mathcal{L}_{N_\sigma + 1} (v), \quad  \forall v\in \mathcal{X}^{\rm hf}_{\rm bc} \quad \text{with} \quad \mathcal{L}_{N_\sigma + 1}^{(k)} = \int_{\Omega} f_v^{(k)}\cdot v \ dx  + \int_{\Gamma_n} f_s^{(k)} \cdot v \ ds $$

This leads to the modification of the Gramian matrix for the last column and the last row:

$$ 
\left\{
\begin{array}{rclc}
&\left(\bm{\Sigma}_N^{(k)}\right)_{n,m}&=\left(\bm{\Sigma}_N\right)_{n,m}, &\forall n,m\in \{1,...,N_\sigma\} \\
&\left(\bm{\Sigma}_N^{(k)}\right)_{n,N_\sigma+1}&=\left(\psi_n^{\sigma}, \ \psi_{N_\sigma + 1}^{\sigma, (k)}\right), &\forall n\in   \{1,...,N_\sigma+1\}
\end{array}
\right.
$$

In practice, we can observe that the $N_\sigma\times N_\sigma$ upper-left submatrix doesn't change over time. A cost-efficient implementation of the Gramian matrix would be only to change the appropriate row over time. Furthermore, we can also observe that for proportionnal loadings (often used for numerical examples in elasto-plasticity, one can compute only one Riesz element and multiply by the appropriate constant at each timestep).

\subsection{Normalisation of the error indicator}

In order not to have values of dual norms that differ depending on the order of magnitude of the loading, we choose to normalize the residual using the norm of the Riesz elements for the external loadings. Moreover, this choice seems consistent with the relative convergence criteria used in pratice in \textsf{code$\_$aster} (see Appendix \ref{app:newton:subsec:stopcrit}). We define $\bm{\widetilde{\Sigma}}_N\in \mathbb{R}^{N_\sigma + 1, N_\sigma + 1}$:

$$\forall n,m\in \{1,...,N_\sigma\}, \quad \left(\bm{\widetilde{\Sigma}}_N^{(k)}\right)_{n,m} = \frac{\left(\bm{\Sigma}_N^{(k)}\right)_{n,m}}{\left(\bm{\Sigma}_N^{(k)}\right)_{N_\sigma+1,N_\sigma+1}} = \frac{\left(\bm{\Sigma}_N^{(k)}\right)_{n,m}}{\norm{\psi_{N_\sigma + 1}^{\sigma, (k)}}^2}$$

The actual error indicator used in our computations is:

\begin{equation}
\Delta_{N, \mu}^{\rm av} = \sqrt{\frac{1}{K} \sum\limits_{k=1}^K \left( \Delta_{N, \mu}^{(k)} \right)^2}, \quad \text{with} \quad \left(\Delta_{N, \mu}^{(k)} \right)^2= \left(\widetilde{\mathbf{\alpha}}^{(k)}_{\sigma, \mu}\right)^T\cdot \bm{\widetilde{\Sigma}}_N^{(k)}\cdot \widetilde{\mathbf{\alpha}}^{(k)}_{\sigma, \mu}
\end{equation}

\section{Details about the elastoplastic solver}
\label{app:elastoplasticity}

The purpose of this appendix is to supply the stages of the numerical procedure adopted so that the work can be reproduced. Plasticity comes to the proficiency of a material to undergo irreversible deformations in reaction to an applied loading. Likewise, elastoplasticity refers to a behaviour where the material has several response regimes: a plastic behaviour for 'small' loadings, and a plastic behaviour (permanent deformations) over some loading amplitude. 

\subsection{Incremental algorithm}

We provide here the choice of the time discretization algorithm used to solve the physical problem detailed in the section \ref{sec:model}. We rely on the discretisation schemes presented by the Eq.\eqref{sec:intro:subsec:context:mechanicalPb:quasistatic}. The time integration of the mechanical behavior of the problem is performed from the computation of a deformation increment:

\begin{equation}
\varepsilon^{(k)} = \varepsilon^{(k-1)} + \Delta \varepsilon^{(k-1)}
\end{equation}

We recall (see Eq. \eqref{sec:Phys:subsec:Plas:eq:constitutiveElas:notations}) that $e$ (resp. $s$) stands for the deviator of the strain (resp. stress) tensor. The discretization of the problem boils down to finding $\left(\Delta p^{(k-1)}, \ \Delta \varepsilon^{p, (k-1)}\right)$ for a given $\Delta \varepsilon^{(k-1)}$  such that:
\begin{equation*}
\left\{
\begin{array}{rcl}
p^{(k)} &=& p^{(k-1)} + \Delta p^{(k-1)}\\
\varepsilon^{p, (k)} &=& \varepsilon^{p, (k-1)} + \Delta \varepsilon^{p, (k-1)}
\end{array}
\right.
\quad \text{with,}
\end{equation*}
\begin{equation}
\left\{
\begin{array}{rcl}
\sigma^{(k)} &=& \displaystyle \sigma^{(k-1)} + \frac{E \nu}{\left(1+\nu\right) \left(1-2\nu\right)}\text{Tr}\left(\Delta \varepsilon^{(k)}\right) + \frac{E}{1+\nu}\left(\Delta e^{(k)} - \Delta \varepsilon^{p, (k)} \right)\\ \\
\sigma^{{\rm eq}, (k)}& - & R\left(p^{(k-1)} + \Delta p^{(k-1)}\right) \leq 0 \\ \\
\Delta \varepsilon^{p, (k-1)} &=& \displaystyle \Delta p^{(k-1)}\frac{3}{2\sigma^{{\rm eq}, (k)}} s^{(k)} \\ \\
\Delta p^{(k)} &\geq & 0 \\ \\
\Delta p^{(k)}& &\left[\sigma^{{\rm eq}, (k)} - R\left(p^{(k-1)} + \Delta p^{(k-1)} \right)\right] = 0
\end{array}
\right.
\end{equation}

We choose to consider an algorithm referred to as incremental in the literature, with a first-order accurate time discretization. The solution varies depending on whether the evolution is exclusively elastic or elastoplastic. Such a procedure adopted is referred to as the return mapping algorithm (or radial return)\cite{wilkins1963calculation}. It resorts to an elastic prediction phase, where the stress field is derived under the assumption of a purely elastic material ($\sigma^{\rm elas}_{n+1}$). The function $f\left(\sigma^{\rm elas}_{n+1}, \ p_n\right)$ is then estimated based on this prediction. If the solution obtained remains in the elastic region, the next iteration can be launched. Otherwise, a correction is performed by solving the nonlinear equation:

\begin{equation}
\label{sec:Phys:subsec:Incr:eq:consistency}
f\left(\sigma_{n+1}, p_n + \Delta p_n\right) = 0
\end{equation}

This equation is nonlinear and is solved through a Newton solver (secant method). The set of unknowns is inferred from the plastic deformation increment $\Delta p_n$. Note that this algorithm is even applied for static problems. In this case, a pseudo-time is introduced. From a physical perspective, it can be understood as a time modeling the evolution of the irreversibility within the material.

\begin{algorithm}[h!]
\caption{Return mapping algorithm}\label{alg:cap}
\begin{algorithmic}
\State Computation of the elastic prediction $s_{n+1}^{\rm elas}= s_{n} + 2\mu\Delta e_n$
\State Stress computation  $\sigma^{\rm elas}_{n+1}$, $\sigma^{\rm elas, \rm eq}_{n+1}$
\State Computation of the criterium $f\left(\sigma^{\rm elas}_{n+1}, \ p_n\right)$
\If{$f\left(\sigma^{\rm elas}_{n+1}, \ p_n\right) \leq 0$} \Comment{ Elastic Evolution}
    \State Computation of the stress and internal variables:
    $$\sigma_{n+1}=\sigma_{n+1}^{\rm elas}, \quad \varepsilon_{n+1}^p=\varepsilon_{n}^p, \quad p_{n+1}=p_n $$
\EndIf
\If{$f\left(\sigma^{\rm elas}_{n+1}, \ p_n\right) > 0$} \Comment{Elastoplastic}
    \State Find $\Delta p_n$ solution of\Comment{Eq.\eqref{sec:Phys:subsec:Incr:eq:consistency}}
    $$\sigma^{\rm elas, \rm eq}_{n+1} - \frac{3 E}{2(1+\nu)} \Delta p_n - R\left(p_n + \Delta p_n\right) = 0$$
    \State Computation of plastic deformation increments:
    $$\varepsilon_{n+1}^p =\varepsilon_{n+}^p+\Delta \varepsilon_{n}^p, \qquad p_{n+1} = p_n + \Delta p_n  $$
    \State Update of stress and internal variable:\Comment{Eq.\eqref{sec:Phys:subsec:Plas:eq:constitutiveElas}}
    $$\sigma_{n+1}= \sigma_{n} +  \frac{E\nu}{\left(1+\nu\right) \left(1-2\nu\right)}\text{Tr}\left(\Delta \varepsilon_n\right)\mathds{1} + \frac{E}{1+\nu}\left(\Delta e_n - \Delta e_n^p\right)$$
\EndIf
\end{algorithmic}
\end{algorithm}

\end{document}